\documentclass[reqno,draft, 12pt]{amsart}
\usepackage{amsmath}
\usepackage{amsfonts}
\usepackage{enumerate}
\usepackage{amssymb}
\usepackage[all]{xypic}

\DeclareFontEncoding{OT2}{}{}

\numberwithin{equation}{section}

\newcommand{\rar}{\longrightarrow}
\newcommand{\rarab}[1]{\overset{#1}{\longrightarrow}}

\newcommand{\al}{\alpha}
\newcommand{\be}{\beta}
\newcommand{\ga}{\gamma}
\newcommand{\Ga}{\Gamma}
\newcommand{\de}{\delta}

\newcommand{\la}{\lambda}
\newcommand{\La}{\Lambda}

\newcommand{\ka}{\kappa}
\newcommand{\eps}{\epsilon}
\newcommand{\sg}{\sigma}

\newcommand{\om}{\omega}

\newcommand{\vp}{\varpi}

\newcommand{\bC}{{\mathbb C}}

\newcommand{\bN}{{\mathbb N}}

\newcommand{\bQ}{{\mathbb Q}}
\newcommand{\bR}{{\mathbb R}}

\newcommand{\bZ}{{\mathbb Z}}

\newcommand{\cO}{{\mathcal O}}

\newcommand{\cW}{{\mathcal W}}

\newcommand{\limproj}{\lim\limits_{\longleftarrow}}

\newcommand{\Hom}{\operatorname{Hom}}

\newcommand{\Spec}{\operatorname{Spec}}

\newcommand{\id}{\operatorname{id}}
\newcommand{\im}{\operatorname{im}}

\newcommand{\rank}{\operatorname{rank}}

\newcommand{\Gal}{\operatorname{Gal}}
\newcommand{\GL}{\operatorname{GL}}
\newcommand{\spin}{\operatorname{sp}}
\newcommand{\tens}{\otimes}

\newcommand{\indlim}[1]{\lim\limits_{\underset{N}{\longrightarrow}}}
\newtheorem{thm}{Theorem}[section]
\newtheorem{cor}[thm]{Corollary}
\newtheorem{lem}[thm]{Lemma}
\newtheorem{sublem}[thm]{Sublemma}
\newtheorem{prop}[thm]{Proposition}

\newtheorem{rem}[thm]{Remark}

\newcommand{\Kb}{\overline{K}}

\newcommand{\Fbv}{\overline{F}_v}
\newcommand{\gk}{\Gal(\Kb/K)}

\newcommand{\gkl}{\Gal(\Kb/L)}
\newcommand{\wk}{\cW(\Kb/K)}
\newcommand{\wdk}{\cW'(\Kb/K)}
\newcommand{\wdl}{\cW'(\Fbv/F_v)}
\newcommand{\wel}{\cW(\Fbv/F_v)}
\newcommand{\gku}{\Gal(\Kb/K^{unr})}
\newcommand{\gks}{\Gal(\overline{k}/k)}
\newcommand{\wf}{\cW(\Kb/F)}
\newcommand{\wl}{\cW(\Kb/L)}
\newcommand{\wkl}{\cW(L/K)}

\newcommand{\ind}{\text{Ind}}
\newcommand{\res}{\text{Res}}

\newcommand{\ktt}{(K^{\times})^r}
\newcommand{\kbtt}{(\Kb^{\times})^r}
\newcommand{\kbt}{\Kb^{\times}}
\newcommand{\ot}{\cO^{\times}}
\newcommand{\ott}{(\cO^{\times})^r}
\newcommand{\kt}{K^{\times}}
\newcommand{\bbe}{\begin{equation}}
\newcommand{\ee}{\end{equation}}
\newcommand{\ol}{(\cO^{\times})^r\cap\La}

\newcommand{\pair}[1]{\langle #1\rangle}

\title[Root numbers of curves]{Root numbers of curves}
\author[M.~Sabitova]{Maria Sabitova\address{Department of Mathematics,
\hfill\newline University of Pennsylvania, \hfill\newline Philadelphia,
PA 19104, \hfill\newline e-mail: sabitova@math.upenn.edu}}

\date{October 30, 2004}

\begin{document}

\begin{abstract}
We generalize a theorem of D. Rohrlich concerning root numbers of elliptic curves over the field of rational numbers. Our result applies to curves of all higher genera over number fields. Namely, under certain conditions which naturally extend the conditions used by D. Rohrlich, we show that the root number $W(X,\tau)$ associated to a smooth projective curve $X$ over a number field $F$ and a complex finite-dimensional irreducible representation $\tau$ of $\Gal(\overline{F}/F)$ with real-valued character is equal to 1. In the case where the ground field is $\bQ$, we show that our result is consistent with the refined version of the conjecture of Birch and Swinnerton-Dyer.
\end{abstract}

\maketitle

\setcounter{tocdepth}{1}

\section{Introduction}\label{s:intro}
The main object of study in our paper is the root number
$W(X,\tau)$ associated to a smooth projective curve $X$ of genus $g$ over a
number field $F$ and a continuous irreducible complex finite-dimensional
representation $\tau$ of $\Gal(\overline{F}/F)$ with real-valued
character. The root number $W(X,\tau)$ is a complex number of
absolute value $1$. Assume for simplicity that $F=\bQ$. Then
$W(X,\tau)$ appears in the following conjectural functional
equation:
\bbe\label{eq:fun}
\La(J_X,\tau,s)=W(X,\tau)\cdot\La(J_X,\tau^*,2-s),
\ee
where $s\in\bC$, $J_X$ is the Jacobian of $X$, $\tau^*$ is the dual to
$\tau$, and $$\La(J_X,\tau,s)=A_{\tau}^s\cdot\Ga(s)^{g\dim\tau}\cdot
L(J_X,\tau,s)$$ for some positive constant $A_{\tau}$ and the twisted
$L$-function $L(J_X,\tau,s)$ which is a meromorphic function of $s$ defined in a right half-plane. This function is conjectured to have an analytic continuation to the entire complex plane. Since $\tau$ has
real-valued character, $\tau\cong\tau^*$ and $W(X,\tau)=\pm
1$. Moreover, assuming \eqref{eq:fun} and considering the
power series expansion of $L(J_X,\tau,s)$ about $s=1$, we get:
\bbe\label{eq:k}
W(X,\tau)=(-1)^{\text{ord}_{s=1}L(J_X,\tau,s)}.
\ee

In our paper we generalize one result by D. Rohrlich for elliptic
curves (\cite{r1}, p. 313, Prop. E) to smooth projective
curves of an arbitrary genus. We prove the following theorem:

\begin{thm}\label{th:main0}
Let $F$ be a number field, $L$ a finite Galois extension of $F$,
and $\tau$ an irreducible complex finite-dimensional
representation of $\Gal(L/F)$ with real-valued character. Let $g$
be a fixed positive integer and assume that the decomposition
subgroups of $\Gal(L/F)$ at all the places of $F$ lying over all
the primes less or equal to $2g+1$ are abelian. If the Schur index
$m_{\bQ}(\tau)$ is $2$ 
then $W(X,\tau)=1$ for
every smooth projective curve $X$ of genus $g$ over $F$.
\end{thm}

If $F=\bQ$ then Theorem \ref{th:main0} can be predicted by the
conjectures of Birch-Swinnerton-Dyer and Deligne-Gross. Namely,
the conjectures of Birch-Swinnerton-Dyer and Deligne-Gross imply
\bbe\label{eq:l}
\text{ord}_{s=1}L(J_X,\tau,s)=\langle\sg_X,\tau \rangle,
\ee
where $\sg_X$ is the natural representation of
$\Gal(\overline{\bQ}/\bQ)$ on $\bC\otimes_{\bZ}J_X(\overline\bQ)$ and $\langle\sg_X,\tau \rangle$ is the multiplicity of $\tau$ in $\sg_X$ (\cite{r3}, p. 127, Prop. 2). This equality is sometimes called the refined version of Birch and Swinnerton-Dyer conjecture. Thus, we get from
\eqref{eq:k} and \eqref{eq:l}:
\bbe\label{eq:m}
W(X,\tau)=(-1)^{\langle\sg_X,\tau \rangle}.
\ee
Since $\sg_X$ is realizable over $\bQ$ and $\tau$ is irreducible,
$m_{\bQ}(\tau)$ 
divides $\langle\sg_X,\tau \rangle$. Thus, if
$m_{\bQ}(\tau)=2$ 
then $W(X,\tau)=1$ for every smooth projective curve $X$ over $\bQ$ (cf. \cite{r1}, p. 313).

To prove Theorem \ref{th:main0} we use the following formula:
\bbe\label{eq:global root}
W(X,\tau)=\prod_v W(X_v,\tau_v),
\ee
where $v$ runs through all the places of $F$, $X_v=X\times_F F_v$,
$F_v$ denotes the completion of $F$ with respect to $v$, and
$\tau_v$ is the restriction of $\tau$ to
$\Gal(\overline{F_v}/F_v)\hookrightarrow\Gal(\overline F/F)$. To
define $W(X_v,\tau_v)$ for every place $v$ let $\sg'_v$ denote the
representation of the Weil-Deligne group $\wdl$ associated to
$X_v$. Then $W(X_v,\tau_v)=W(\sg'_v\otimes\tau_v)$, where $\tau_v$
is viewed as a representation of $\wdl$. We will in fact show the following stronger result:
\begin{thm}\label{th:main}
$W(X_v,\tau_v)=1$ for all $v$ under the hypotheses of Theorem \ref{th:main0}.
\end{thm}

First, we describe $W(X_v,\tau_v)$. If $v$ is an infinite place then
\bbe\label{eq:n}
W(X_v,\tau_v)=(-1)^{g\dim\tau}
\ee
(Lemma \ref{l:archimed}).

If $v$ is a finite place, then 
\bbe\label{eq:eps}
W(\sg_v'\otimes\tau_v)=\frac{\epsilon(\sg_v'\otimes\tau_v,\psi_v,dx_v)}
{|\epsilon(\sg_v'\otimes\tau_v,\psi_v,dx_v)|},
\ee
where $\psi_v$ is an additive character of $F_v$ and $dx_v$ is a Haar measure on $F_v$ (\cite{r2}, p. 144). Here $\sg'_v$ is isomorphic to the representation of
$\wdl$ afforded by the $l$-adic cohomology groups of $X_v$
(\cite{chin}, p. 427), where $l$ is a rational prime, $l\ne\text{char}(k_v)$, and $k_v$ is the residue class field of $F_v$. Since $X_v$ is a smooth projective
curve, $\sg_v'$ is the representation of $\wdl$ associated to the
natural $l$-adic representation of $\Gal(\Fbv/F_v)$ on
$H^1_l(X_v)$, which is the $\bQ_l$-dual to
$V_l(J_v)=T_l(J_v)\otimes_{\bZ_l}\bQ_l$, where $J_v$ is the
Jacobian of $X_v$, and $T_l(J_v)$ is the $l$-adic Tate module of
$J_v$. Clearly, $W(\sg_v'\otimes\tau_v)$ does not depend on the choice of $dx_v$ and it turns out that $W(\sg_v'\otimes\tau_v)$ does not depend on the choice of $\psi_v$ either. Moreover, $W(\sg_v'\otimes\tau_v)=\pm 1$ (see Subsection \ref{sub:loc gen}). 

We consider two cases: $J_v$ is an abelian variety with
potential good reduction, and the general case. If $J_v$ has
potential good reduction, it follows from the theory of
Serre-Tate that $\sg_v'$ is actually a representation of the Weil group $\wel$
(see Subsection \ref{sus: potential}). If
$\text{char}(k_v)>2g+1$, we use the theory of Serre-Tate together
with methods of the representation theory to describe the class of
$\sg'_v\otimes\om_v^{1/2}$ in the Grothendieck group of virtual
representations of $\wel$ (Corollary \ref{cor:ad}, Formula
\eqref{eq:realiz2}). Here $\om_v$ is the one-dimensional
representation of $\wel$ given by
$$
\omega_v\vert_{I_v}=1,\quad \omega_v(\Phi_v)=q_v^{-1},
$$
where $I_v$ is the inertia group of $\Gal(\Fbv/F_v)$, $\Phi_v$ is
an inverse Frobenius element of $\Gal(\Fbv/F_v)$, and
$q_v=\text{card}(k_v)$. Since the root number of representations
of $\wel$ is multiplicative in short exact sequences, this result
enables us to prove the following formula for
$W(\sg_v'\otimes\tau_v)$ when $\text{char}(k_v)>2g+1$ (Proposition
\ref{pr:gen of Rohrlich}):
\bbe\label{eq:o}
W(\sg_v'\otimes\tau_v)=\det\mu(-1)^{\dim\tau}\cdot\det\tau(-1)^{l_1}\cdot\al^{\dim\tau}\cdot(-1)^{l_2},
\ee
where $l_1=\dim\mu+\frac{1}{2}(\dim\mu_1+\dotsb+\dim\mu_a)$,
$\mu,\mu_1,\ldots,\mu_a$ are some representations of $\wel$,
$\al=\pm 1$, $l_2=a\cdot\langle 1,\tau_v\rangle+
a\cdot\langle\eta_v,\tau_v\rangle+a\cdot\langle
\hat\mu_1\oplus\dotsb\oplus\hat\mu_a,\tau_v\rangle$,
$\hat\mu_1\oplus\dotsb\oplus\hat\mu_a$ is a representation of
$\Gal(\Fbv/F_v)$ realizable over $\bQ$, and $\eta_v$ is the
unramified quadratic character of $F_v^{\times}$ (cf. \cite{r1},
p. 318, Thm. 1).

In  the general case we use the theory of uniformization of
abelian varieties. According to this theory there exists a
semi-abelian variety $G_v$ over $F_v$ and a discrete subgroup
$Y_v$ of $G_v$ such that, in terms of rigid geometry, $J_v$ is
isomorphic to the quotient $G_v/Y_v$. The semi-abelian variety
$G_v$ fits into an exact sequence
\bbe\label{eq:kl}
0\longrightarrow T_v\longrightarrow G_v\rarab{f_v} A_v
\longrightarrow 0,
\ee
where $A_v$ is an abelian variety over $F_v$ with potential good
reduction, $T_v$ is a torus over $F_v$ of dimension $r$; $Y_v$ is
an \'etale sheaf of free abelian groups over $\Spec(F_v)$ of rank
$r$. To describe $\sg_v'$ in this case we use a formula
of M. Raynaud (\cite{ray}, p. 314) which gives the action of the inertia group $I_v$ on the $l^n$-torsion points of an abelian variety over a local field in the case when the uniformization data splits. We need this formula to show that in this
case
\bbe\label{eq:po}
\sg'_v\cong\ka_v\oplus(\chi_v\otimes\om_v^{-1}\otimes\spin(2)),
\ee
where $\ka_v$ is the representation of $\wdl$ associated to the
natural $l$-adic representation of $\Gal(\Fbv/F_v)$ on the dual
space to $V_l(A_v)=T_l(A_v)\otimes_{\bZ_l}\bQ_l$,
$$
\chi_v:\Gal(\Fbv/F_v)\rar\GL_r(\bZ)
$$
is the representation of $\Gal(\Fbv/F_v)$ corresponding to the
Galois module $Y_v(\overline F_v)$, and $\spin(2)$ is given by \eqref{eq:special} (see Proposition \ref{pr:exact_sequence}). Since the root number of a direct sum of
representations of $\wdl$ equals the product of the root numbers
of the summands, we get from \eqref{eq:po}
\bbe\label{eq:r}
W(\sg'_v\otimes\tau_v)=W(\ka_v\otimes\tau_v)\cdot
W(\chi_v\otimes\om_v^{-1}\otimes\tau_v\otimes\spin(2)).
\ee
We show that if $\text{char}(k_v)>2g+1$ then \eqref{eq:o} holds
for $\ka_v$, i.e.
\bbe\label{eq:j}
W(\ka_v\otimes\tau_v)=\det\mu(-1)^{\dim\tau}\cdot\det\tau(-1)^{l_1}\cdot\al^{\dim\tau}\cdot(-1)^{l_2},
\ee
where $l_1=\dim\mu+\frac{1}{2}(\dim\mu_1+\dotsb+\dim\mu_a)$,
$\mu,\mu_1,\ldots,\mu_a$ are some representations of $\wel$,
$\al=\pm 1$, $l_2=a\cdot\langle 1,\tau_v\rangle+
a\cdot\langle\eta_v,\tau_v\rangle+a\cdot\langle
\hat\mu_1\oplus\dotsb\oplus\hat\mu_a,\tau_v\rangle$,
$\hat\mu_1\oplus\dotsb\oplus\hat\mu_a$ is a representation of
$\Gal(\Fbv/F_v)$ realizable over $\bQ$, and $\eta_v$ is the
unramified quadratic character of $F_v^{\times}$.

The rest of the proof of Theorem \ref{th:main} is analogous to one of
Proposition E (\cite{r1}, p. 347). Namely, it follows from Lemma
on p. 339 and Lemma on p. 347 in \cite{r1} that $\dim\tau$ is
even. Hence we get from \eqref{eq:n} that $W(X_v,\tau_v)=1$ for
infinite places. If $v$ is a finite place then the assumption $m_{\bQ}(\tau)=2$
implies that $W(\chi_v\otimes\om_v^{-1}\otimes\tau_v\otimes\spin(2))=1$ (\cite{r1}, p. 327, Prop. 6), hence we have from \eqref{eq:r}
\bbe\label{eq:koko}
W(\sg'_v\otimes\tau_v)=W(\ka_v\otimes\tau_v).
\ee

If $v$ is a finite place such that
$\text{char}(k_v)>2g+1$ then \eqref{eq:j} holds which, together with the assumption $m_{\bQ}(\tau)=2$, implies that $W(\ka_v\otimes\tau_v)=1$. Hence
$W(\sg'_v\otimes\tau_v)=1$ by \eqref{eq:koko}.

If $v$ is a finite place such that
$\text{char}(k_v)\leq 2g+1$ then the conditions on bad primes in
Theorem \ref{th:main0} imply that $\tau_v$ is symplectic. We show
that $\ka_v\otimes\om_v^{1/2}$ is symplectic too (Corollary
\ref{c:kasympl}). Since the real powers of $\om_v$ do not change
the root number,
\bbe\label{eq:lolo}
W(\ka_v\otimes\tau_v)=W(\ka_v\otimes\om_v^{1/2}\otimes\tau_v),
\ee
where $W(\ka_v\otimes\om_v^{1/2}\otimes\tau_v)=1$ as the root number of the tensor product of two symplectic representations of $\wel$ (\cite{r1}, p. 319, Prop. 2 and the remark after it). Thus,
in this case we also have $W(X_v,\tau_v)=1$ by \eqref{eq:koko} and \eqref{eq:lolo}.

\smallskip
This paper is organized in the following way. In Section
\ref{s:local} we study the root number $W(\sg'\otimes\tau)$, where
$\tau$ is a complex finite-dimensional representation of $\gk$
with real-valued character, $K$ is a local non-Archimedean field
of characteristic zero, and $\sg'$ is the representation of $\wdk$
associated to the $l$-adic representation of $\gk$ on the dual
$V_l(M)^*$ to $V_l(M)$, where $M$ is an abelian variety over $K$.
Subsection \ref{sub:loc gen} contains general facts and
notation used in the paper. In Subsection \ref{sus: potential} we
study the case of an abelian variety with potential good
reduction. Subsection \ref{sec:gen}  deals with the general case.
In Section \ref{sec:th} we give the proof of Theorem
\ref{th:main}. We put proofs of the results of Subsection
\ref{sus: potential} in Appendix \ref{sec:b}. Appendix \ref{ap:c}
contains two lemmas needed for the proof of the main result of
Subsection \ref{sec:gen} (Proposition \ref{pr:exact_sequence}). In
Appendix \ref{ap:a} we give a description of the representation of
$\wdk$ associated to the natural $l$-adic representation of $\gk$
on $V_l(M)^*$ in the case when $M$ is the quotient of a torus by a
discrete subgroup (Proposition \ref{p:nonsplit tor}). This result
is a special case of Proposition \ref{pr:exact_sequence}, however
we give an elementary proof which does not rely on the result by
M. Raynaud mentioned above. In fact, we prove a more general result, which can
be used to prove formula \eqref{eq:po} for an arbitrary Jacobian variety $J_v$
over $F_v$ without using Raynaud's
result in a special case when in 
\eqref{eq:kl} the image of $Y_v$
under $f_v$ is finite. Instead of Raynaud's
result the description of symplectic admissible representations of $\wdl$ can be used in this case. We give a
description of unitary, orthogonal, or symplectic admissible
representations of $\wdk$ in Appendix \ref{ap:d}.

Unless stated otherwise, we assume that all the representations
under consideration are complex and finite-dimensional.
\medskip

\noindent
{\bf Acknowledgments}. I would like to thank my adviser Ted
Chinburg for suggesting the problem and useful discussions. I am
also grateful to Siegfried Bosch, Ching-Li Chai, Robert Kottwitz,
and Michel Raynaud for answering my questions related to the
uniformization theory, and especially to Michel Raynaud for
providing the reference \cite{ray}.

\section{Local case}\label{s:local}
\subsection{General facts and notation}\label{sub:loc gen}

Let $K$ be a non-Archimedean local field of characteristic zero
with ring of integers $\cO$, residue class field $k$, and a
uniformizer $\varpi$. Let $\Kb$ be a fixed algebraic closure of
$K$ and let $K^{unr}$ be the maximal unramified extension of $K$
contained in $\Kb$. Let $I=\gku$ be the inertia subgroup of $\gk$
and let $\Phi$ be an inverse Frobenius element of $\gk$, i.e. $\Phi$ is a
preimage of the inverse of the Frobenius automorphism under the
decomposition map

$$
\pi:\gk \longrightarrow \gks.
$$

By a representation $\sg$ of the Weil group $\wk$ we mean a continuous homomorphism
$$
\sg:\wk\rar\GL(U),
$$
where $U$ is a finite-dimensional complex vector space.
Let $\omega:\wk\longrightarrow \bC^{\times}$ be the
one-dimensional representation of $\wk$ given by
$$
\omega\vert_{I}=1,\quad \omega(\Phi)=q^{-1},
$$
where $q=\text{card}(k)$. For a finite extension $F$ of $K$
contained in $\Kb$, we identify by local class field theory the
one-dimensional representations of $\wf$ with characters of
$F^{\times}$ (i.e. continuous homomorphisms from $F^{\times}$ into
$\bC^{\times}$). Also, if $\phi$ is a representation of $\wf$, the
representation of $\wk$ induced by $\phi$ will be denoted by
$\ind^F_K\phi$. Analogously, if $\psi$ is a representation of
$\wk$, then the restriction of $\psi$ to $\wf$ will be denoted by
$\res^F_K\psi$.

By a representation $\sg'$ of the Weil-Deligne group $\wdk$ we mean a continuous homomorphism
$$
\sg':\wdk\rar\GL(U),
$$
where $U$ is a finite-dimensional complex vector space and the restriction of $\sg'$ to the subgroup $\bC$ of $\wdk$ is complex analytic. It is known that there is a bijection between representations of $\wdk$ and pairs $(\sg,N)$, where $\sg:\wk\rar\GL(U)$ is a representation of $\wk$ and $N$ is a nilpotent endomorphism on $U$ such that
$$
\sg(g)N\sg(g)^{-1}=\om(g)N,\quad g\in\wk.
$$
In what follows we identify $\sg'$ with the corresponding pair $(\sg,N)$ and write $\sg'=(\sg,N)$. Also, a representation $\sg$ of $\wk$ is identified with the representation $(\sg,0)$ of $\wdk$ (\cite{r2}, p. 128, \S 3). 

\medskip

For a positive integer $n$ let $\spin(n)=(\sg,N)$ denote the special representation of dimension $n$, i.e. the representation of $\wdk$ on $\bC^n$ (with the standard basis $e_0,\ldots,e_{n-1}$) given by the following formulas:
\begin{eqnarray}
\sg(g)e_i & = & \om(g)^ie_i, \quad 0\leq i\leq n-1,\,g\in\wk,\label{eq:special}\\
Ne_j & = & e_{j+1},\,\qquad 0\leq j\leq n-2,\nonumber\\
Ne_{n-1} & = & 0. \nonumber
\end{eqnarray}
We say that a representation $\sg'=(\sg,N)$ of $\wdk$ is {\it admissible} if $\sg$ is semisimple (\cite{r2}, p. 132, \S 5).

Let $M$ be an abelian variety over $K$. For a rational prime $l$ different
from $p=\text{char}(k)$ let $T_l(M)$ be the $l$-adic Tate module
of $M$. It is a free $\bZ _l$-module of rank $2g$, where $g=\dim
M$. Put $V_l(M)=T_l(M)\otimes_{\bZ _l} \bQ_l$ and denote by
$V_l(M)^*$ the dual of $V_l(M)$. Let
$$
\sigma_l:\gk\longrightarrow \GL(V_l(M)^*)
$$
denote the contragredient of the natural $l$-adic representation
of $\gk$ on $V_l(M)$. We are interested in the representation
$\sigma'=(\sigma,N)$ of $\wdk$ associated
to $\sigma_l$. Let $\imath:\bQ_l\hookrightarrow \bC$ be a field
embedding. Then $\sigma:\wk\longrightarrow
\GL(V_l(M)^*\otimes_{\imath}\bC)$ is a representation of $\wk$
(which is not necessarily obtained from the restriction of $\sg_l$ to $\wk$ by extending scalars via $\imath:\bQ_l\hookrightarrow\bC$) and
$N\in\text{End}(V_l(M)^*\otimes_{\imath}\bC)$ is a nilpotent
endomorphism (see \cite{r2}, p. 130, \S4 for more detail).
A priori, $\sg'$ depends on the choice of $l$ and $\imath$, but by abuse of notation we write $\sg'$ instead of $\sg'_{l,\imath}$. We will prove later that in our context $\sg'$ does not depend on the choice of $l$ and $\imath$.
Let $\tau$ be a representation of $\gk$ with real-valued
character. Our goal in this section is to compute the root number
$W(\sg'\otimes\tau)$. We are particularly concerned with the case
when $M$ is the Jacobian of a smooth projective curve over $K$. In
this case $M$ is self-dual and, consequently, there is the Weil
pairing on $M_{l^n}=\text{Hom}(\bZ/l^n\bZ,M(\Kb))$, which induces
a nondegenerate, symplectic, $\gk$-equivariant pairing
$$
\langle -,-\rangle:T_l(M)\times T_l(M)\longrightarrow
\bZ_l\otimes\omega_l,
$$
where $\omega_l$ is the $l$-adic cyclotomic character of $\gk$.
It is easy to show that $\sg'\otimes\omega^{1/2}$ is
symplectic (cf. \cite{r2}, p. 150, \S 16). This is the only use we
will make of the fact that $M$ is the Jacobian variety of a curve.
Thus, throughout this section we assume that $M$ is an abelian
variety such that $\sg'\otimes\omega^{1/2}$ is symplectic. Then $\sg'\otimes\omega^{1/2}\otimes\tau$ is self-contragredient and of trivial determinant, hence $W(\sg'\otimes\tau)=W(\sg'\otimes\omega^{1/2}\otimes\tau)$, where $W(\sg'\otimes\omega^{1/2}\otimes\tau)$ does not depend on the choice of an additive character of $K$  and $W(\sg'\otimes\omega^{1/2}\otimes\tau)=\pm 1$ (\cite{r1}, p. 315).

\smallskip

The main theory we are using to find a formula for
$W(\sg'\otimes\tau)$ is the theory of uniformization of abelian
varieties. According to this theory there exists a semi-abelian
variety $G$ over $K$ and a discrete subgroup $Y$ of $G$ such that,
in terms of rigid geometry, $M$ is isomorphic to the quotient
$G/Y$. The semi-abelian variety $G$ fits into an exact sequence
\bbe\label{eq:uniform}
0\longrightarrow T\longrightarrow G\rarab{f} A \longrightarrow 0,
\ee
where $A$ is an abelian variety over $K$ with potential good
reduction, $T$ is a torus over $K$ of dimension $r$; $Y$ is an
\'etale sheaf of free abelian groups over $\Spec(K)$ of rank $r$.

\subsection{Case of an abelian variety with potential good reduction.}\label{sus: potential}
Let $A$ be an abelian variety over $K$ with potential good
reduction and let
$$
\ka_l:\gk\longrightarrow \GL(V_l(A)^*)
$$
denote the natural $l$-adic representation of $\gk$ on $V_l(A)^*$.
First, let us show that the representation $\ka'=(\ka,S)$ of
$\wdk$ associated to $\ka_l$ is actually a representation of
$\wk$, i.e. $S=0$. Indeed, $\ka'$ is a representation of $\wk$ if and only if $\ka_l$ is trivial on an open subgroup of $I$ (\cite{r2}, p. 131, Prop.(i)). Let
$$
\psi_l:\gk\longrightarrow \text{Aut}(T_l(A))
$$
denote the representation corresponding to the $\gk$-module
$T_l(A)$. Since $A$ has potential good reduction, the image by
$\psi_l$ of $I$ is finite (\cite{st}, p. 496, Thm. 2(i)), which
implies that the image by $\ka_l$ of $I$ is finite, hence $\ka_l$
is trivial on an open subgroup of $I$ (cf. \cite{r2}, p. 148).

It is known that a complex finite-dimensional representation $\la$
of a group is semisimple if and only if its restriction to a
normal subgroup of finite index is semisimple (\cite{c}, p. 82,
Prop. 1 and \cite{r2}, p. 148). Moreover, since every subgroup of
finite index contains a normal subgroup of finite index, this
implies that $\la$ is semisimple if and only if its restriction to
a subgroup of finite index is semisimple. This fact will be used frequently throughout this paper, so for convenience we state it as a separate lemma:

\begin{lem}\label{l:find}
A complex finite-dimensional representation 
of a group is semisimple if and only if its restriction to a
subgroup of finite index is semisimple.
\end{lem}

\begin{lem}\label{l:semisimple1}
$\ka$ is semisimple.
\end{lem}

\begin{proof}
Since the image by $\ka$
of $I$ is finite, by Lemma \ref{l:find} it is enough to show that
$\ka(\Phi)$ is diagonalizable. Also, if $L\subset \Kb$ is a finite
extension of $K$ over which $A$ acquires good reduction, then by
the above fact $\ka$ is semisimple if and only if its restriction
to $\wl$ is semisimple. Thus, we can assume that $A$ has good
reduction (cf. \cite{r2}, p. 148). Let $A_0$ be the N\'eron
minimal model of $A$ and $\tilde A=A_0\times_{\cO}k$ the special
fiber of $A_0$. Since $A$ has good reduction, the reduction map
defines a $\gk$-equivariant isomorphism of $T_l(A)$ onto
$T_l(\tilde{A})$, where $\gk$ acts on $T_l(\tilde{A})$ via the
decomposition map $\pi$ (\cite{st}, p. 495, Lem. 2). Thus,
\bbe\label{eq:good red}
V_l(A)\cong V_l(\tilde{A})
\ee
as $\gk$-modules. By Tate's result on Tate's conjecture, the
natural $l$-adic representation $\be_l$ of $\gks$ on
$V_l(\tilde{A})$ is semisimple. Since $\gks$ is abelian, $\be_l$
is a direct sum of one-dimensional representations, hence
$\be_l(\pi(\Phi))$ is diagonalizable, consequently,
$\ka_l^*(\Phi)$ is diagonalizable, because $\ka_l^*(\Phi)$ is
equivalent to $\be_l(\pi(\Phi))$ via \eqref{eq:good red}. This
proves that $\ka(\Phi)$ is diagonalizable, because $\ka(\Phi)$ is
just $\ka_l(\Phi)$ considered as an element of
$\GL(V_l(A)^*\otimes_{\imath}\bC)$.
\end{proof}

\begin{cor}\label{c:dep}
The representation $\ka$ does not depend on choice of $l$ and
$\imath$.
\end{cor}
\begin{proof}
\cite{r2}, p. 148 and Lemma \ref{l:semisimple1}.
\end{proof}

\begin{rem}
Throughout this subsection we assume that $\ka\otimes\om^{1/2}$ is
symplectic. We will prove  that this is true in our context, i.e.
for an abelian variety $A$ over $K$ with potential good reduction
fitting into exact sequence \eqref{eq:uniform} corresponding to an
abelian variety $M$ which is the Jacobian of a smooth projective
curve over $K$ $($Corollary $\ref{c:kasympl})$.
\end{rem}

Since $A$ has potential good reduction, there exists a minimal
finite subextension $L/K^{unr}$ of $\Kb/K^{unr}$ over which $A$
acquires good reduction. It is a Galois extension and it is tamely
ramified if $p>2m+1$, where $m=\dim A$. Moreover, $\gkl$ is
contained in the kernel of the representation $\psi_l$ (\cite{st},
p. 497, Cor. 2 and p. 498, Cor. 3). Thus, $\ka$ and, consequently
$\ka\otimes\om^{1/2}$, can be considered as representations of the
group
$$
\wkl=\wk/\gkl\cong\Gal(L/K^{unr})\rtimes(\Phi),
$$
where $(\Phi)$ is the infinite cyclic group generated by $\Phi$
(cf. \cite{r1}, p. 331). Throughout this subsection we assume that
$p>2m+1$. Then, under this assumption $B=\Gal(L/K^{unr})$ is a
finite cyclic group of order not divisible by $p$ and
$\ka\otimes\om^{1/2}$ is a semisimple (by Lemma
\ref{l:semisimple1}), symplectic (by assumption) representation of
the semi-direct product $G=B\rtimes(\Phi)$ of finite and infinite
cyclic groups. Using Corollary on p. 499 in \cite{st}, it is
immediate that $\ka$ has $\bQ$-valued character. Since $\om$ is
trivial on $I$, it follows that $\res ^G_B(\ka\otimes\om^{1/2})$
has $\bQ$-valued character. The following results give a
description of such a representation, i.e. a semisimple symplectic
representation $\la$ of a semi-direct product of a finite cyclic
group $B$ and an infinite cyclic group such that the restriction
of $\la$ to $B$ has $\bQ$-valued character. They will be used
later to generalize a formula for the root number obtained by D.
Rohrlich.

\begin{prop}\label{p:sympl irred}
Let $C=(c)$ be an infinite cyclic group generated by an element
$c$ and let $B=(b)$ be a finite cyclic group of order $n$ generated by an
element $b$. Let  $G=B\rtimes C$ be a semi-direct product, where
$C$ acts on $B$ via $c^{-1}bc=b^k$ for some
$k\in(\bZ/n\bZ)^{\times}$. Denote by $s$ the order of $k$ in
$(\bZ/n\bZ)^{\times}$. Then every irreducible symplectic
representation $\la$ of $G$ factors through the group
$H=G/(c^{2s})$ and as a representation of $H$ it has the following
form
$$
\la=\text{Ind}_{B\rtimes\Ga}^H\phi,
$$
where $\Ga$ is a subgroup of $C/(c^{2s})$ generated by an element $c^x$ and
$\phi$ is a one-dimensional representation of $B\rtimes\Ga$
satisfying the following conditions:
\begin{itemize}
\item $\phi(b)=\xi$ for an n-th root of unity $\xi$ of order
$d$ $(d\not= 1,2)$
\item $x$ is the order of $k$ in $(\bZ/d\bZ)^{\times}$
\item $x$ is even
\item $\phi(c^x)=-1$
\item $1+k^{\frac{x}{2}}\equiv0\,(\text{mod}\,d)$.
\end{itemize}
Conversely, every representation of this form is symplectic and
irreducible.
\end{prop}

In the notation of Proposition \ref{p:sympl irred} let
$\la=\text{Ind}_{B\rtimes\Ga}^H\phi$ be a symplectic irreducible
representation of $G$ and $\theta$ the one-dimensional
representation of $B\rtimes\Ga$ such that $\theta(c^x)=-1$,
$\theta(b)=1$. Let
\bbe\label{eq:hat}
\hat\la=\text{Ind}_{B\rtimes\Ga}^H(\phi\otimes\theta).
\ee
Whereas $\la$ is symplectic, $\hat{\la}$ is realizable over $\bR$. It can be checked using Proposition 39 (\cite{S}, p. 109).

For a group $D$ let $R(D)$ denote the Grothendieck group of the
abelian category of finite-dimensional representations of $D$ over
$\bC$. If $\rho$ is such a representation we denote by $[\rho]$
the corresponding element of $R(D)$.

\begin{prop}\label{p:realiz q}
Let $G=B\rtimes C$ be a semi-direct product as in Proposition $\ref{p:sympl irred}$ and $\la$ a semisimple symplectic
representation of $G$. If $\res^G_B\la$ has $\bQ$-valued
character, then
\bbe\label{eq:realiz1}
[\la]=[\mu]+[\mu^{*}]+2\cdot([\mu_0]-[\mu'_0])+[\mu_1]+\dotsb
+[\mu_a],
\ee
where $\mu$ is a representation of $G$, $\mu_0$ and $\mu'_0$ are
symplectic representations of $G$ with finite images,
$\mu_1,\ldots,\mu_a$ are irreducible symplectic subrepresentations
of $\la$ with finite images, $\hat{\mu}_1,\ldots,\hat{\mu}_a$ are
representations with finite images given by \eqref{eq:hat} such
that $\hat{\mu}_1\oplus\dotsb\oplus\hat{\mu}_a$ is realizable over
$\bQ$.
\end{prop}


\begin{cor}\label{cor:ad}
Let $\ka$ be the representation of $\wk$ corresponding to
$V_l(A)^*$ such that $\la=\ka\otimes\om^{1/2}$ is symplectic. Let $m=\dim A$ and
$p>2m+1$. Then
in $R(\wk)$ we have
\bbe\label{eq:realiz2}
[\la]=[\mu]+[\mu^{*}]+2\cdot([\mu_0]-[\mu'_0])+[\mu_1]+\dotsb
+[\mu_a],
\ee
where $\mu$ is a representation of $\wk$, $\mu_0$ and $\mu'_0$ are
symplectic representations of $\wk$ with finite images,
$\mu_1,\ldots,\mu_a$ are irreducible symplectic subrepresentations
of $\la$ with finite images, $\hat{\mu}_1,\ldots,\hat{\mu}_a$ are
representations with finite images given by \eqref{eq:hat} such
that $\hat{\mu}_1\oplus\dotsb\oplus\hat{\mu}_a$ is realizable over
$\bQ$.
\end{cor}

Let $\tau$ be a representation of $\gk$ with real-valued
character. To compute the root number $W(\ka\otimes\tau)$ we
generalize the following result by D. Rohrlich (\cite{r1}, p. 318,
Thm. 1) :
\begin{thm}\label{th:rohrlich}
Let $K$ be a local non-Archimedean field of characteristic zero.
Let $\tau$ be a representation of $\gk$ with real-valued
character. Then
$$
W(\la\otimes\tau)=\det\tau(-1)\cdot\varphi(u_{H_2/K})^{\dim\tau}\cdot(-1)^{\langle
1,\tau\rangle+ \langle\eta,\tau\rangle+\langle
\hat\la,\tau\rangle}.
$$
\end{thm}
\noindent
Here $\la$ is a two-dimensional irreducible, symplectic
representation of $\gk$ of the form $\la=\ind^{H_2}_K\phi$, where
$H_2$ is the unramified quadratic extension of $K$, $\phi$ is a
tame character of $H_2^{\times}$; $\eta$ is the unramified
quadratic character of $K^{\times}$;
$\hat\la=\ind^{H_2}_K(\phi\otimes\theta)$, where $\theta$ is the
unramified quadratic character of $H_2^{\times}$, and
$\varphi(u_{H_2/K})=\pm 1$ (see \cite{r1}, p. 318 for more
detail).

\smallskip

We prove the following generalization of Theorem
\ref{th:rohrlich}:

\begin{prop}\label{pr:gen of Rohrlich}
Let $\ka$ be the representation of $\wk$ corresponding to
$V_l(A)^*$ such that $\ka\otimes\om^{1/2}$ is symplectic. 
Let $m=\dim A$ and $p>2m+1$.
Let $\tau$ be a representation of $\gk$ with real-valued
character. In the notation of Corollary $\ref{cor:ad}$ we have
\bbe\label{pr:pg}
W(\ka\otimes\tau)=\det\mu(-1)^{\dim\tau}\cdot\det\tau(-1)^{l_1}\cdot\al^{\dim\tau}\cdot(-1)^{l_2},
\ee
where $l_1=\dim\mu+\frac{1}{2}(\dim\mu_1+\dotsb+\dim\mu_a)$,
$\al=\pm 1$, $l_2=a\cdot\langle 1,\tau\rangle+
a\cdot\langle\eta,\tau\rangle+a\cdot\langle
\hat\mu_1\oplus\dotsb\oplus\hat\mu_a,\tau\rangle$, and $\eta$ is
the unramified quadratic character of $K^{\times}$.
\end{prop}

\subsection{General case}\label{sec:gen}
We keep the notation of Subsection \ref{sub:loc gen}. 
Let $\sg'=(\sg,N)$
be the representation of $\wdk$ associated to the natural $l$-adic
representation  of $\gk$ on $V_l(M)^*$, $\ka$ the representation
of $\wdk$ associated to the natural $l$-adic representation  of
$\gk$ on $V_l(A)^*$, and
$$
\chi:\gk\longrightarrow \GL_r(\bZ)
$$
the representation corresponding to the $\gk$-module $Y(\Kb)$. It
is known that there is a finite Galois extension $L\subset\Kb$ of
$K$ such that $\Gal(\Kb/L)$ acts trivially on $Y(\Kb)$, hence $\chi$ has
finite image. Here $\ka$ is actually a representation of $\wk$ (see Subsection \ref{sus: potential}) and we identify $\ka$ with the representation $(\ka,0)$ of $\wdk$. Also, we identify $\chi$ with the representation $(\res_{\wk}\chi,0)$ of $\wdk$.

The main result of this subsection is the following proposition:

\begin{prop}\label{pr:exact_sequence}
\bbe\label{eq:sas}
\sg'\cong\ka\oplus(\chi\otimes\om^{-1}\otimes\spin(2)).
\ee
\end{prop}
\noindent
To prove Proposition \ref{pr:exact_sequence} we will need the following lemmas.

\begin{lem}\label{l:ator}
Let $\la'=(\la,R)$ be the representation of $\wdk$ associated to the natural $l$-adic representation of 
$\gk$ on $V_l(T)^*$. Then $R=0$ and
$$
\la\cong\chi\otimes\om^{-1}.
$$
\end{lem}

\begin{proof}
From the exact $\gk$-equivariant sequence \eqref{eq:uniform} we
get the following exact sequence of $\gk$-modules:
$$
0\longrightarrow T(\Kb)\longrightarrow G(\Kb)\rar
A(\Kb)\longrightarrow 0.
$$
Since $T(\Kb)$ is a divisible group, the last sequence induces an
exact $\gk$-equivariant sequence of $l$-adic Tate modules:
$$
0\longrightarrow T_l(T)\longrightarrow T_l(G)\rar T_l(A)
\longrightarrow 0.
$$
By tensoring the above sequence with $\bQ_l$ over $\bZ_l$ and
taking duals over $\bQ_l$ afterwards, we get the exact sequence of
$\gk$-modules:
\bbe\label{eq:111}
0\longrightarrow V_l(A)^*\longrightarrow V_l(G)^*\rar V_l(T)^*
\longrightarrow 0.
\ee

Let $X$ be the character
group of $T$. Then $T(\Kb)\cong\text{Hom}_{\bZ}(X(\Kb),\kbt)$ as
$\gk$-modules over $\bZ$, hence we have the following sequence of
isomorphisms of $\gk$-modules:
\begin{eqnarray}\label{eq:ostar1}
V_l(T) & = & T_l(T)\otimes_{\bZ _l}\bQ _l\cong\text{Hom}_{\bZ}(X(\Kb),T_l(\kbt))\otimes_{\bZ _l}\bQ _l \\
& \cong & \text{Hom}_{\bZ _l}(X(\Kb)\otimes_{\bZ}\bZ
_l,T_l(\kbt))\otimes_{\bZ _l}\bQ _l\nonumber \\ & \cong &
\text{Hom}_{\bQ _l}(X(\Kb)\otimes_{\bZ}\bQ _l,V_l(\kbt))\nonumber
\\ & \cong & (X(\Kb)\otimes_{\bZ}\bQ _l)^*\otimes_{\bQ
_l}V_l(\kbt) \nonumber .
\end{eqnarray}
It is known that there is an injective homomorphism
$\phi:Y\longrightarrow X$ with finite cokernel (\cite{ch}, p. 58),
consequently
$$
Y(\Kb)\otimes_{\bZ}\bQ _l \cong X(\Kb)\otimes_{\bZ}\bQ _l
$$
as $\gk$-modules over $\bQ _l$. Thus, we get from
\eqref{eq:ostar1}
\bbe\label{eq:tensor}
V_l(T)^*\cong(Y(\Kb)\otimes_{\bZ}\bQ _l)\otimes_{\bQ _l}V_l(\kbt)^*.
\ee
Let $\imath:\bQ_l\hookrightarrow \bC$ be a field embedding and let $F_{l,\imath}$ be the functor which associates to an $l$-adic representation of $\gk$ a representation of $\wdk$. Since the image of the representation of $\gk$ on $Y(\Kb)\otimes_{\bZ}\bQ _l$ under $F_{l,\imath}$ is $\chi$ and the image of the representation of $\gk$ on $V_l(\kbt)^*$ under $F_{l,\imath}$ is $\om^{-1}$,
by \eqref{eq:tensor} the image $\la'$ of the representation of $\gk$ on $V_l(T)^*$ under $F_{l,\imath}$ is $\chi\otimes\om^{-1}$. Thus $\la'\cong\chi\otimes\om^{-1}$, hence $\la'$ is a representation of $\wk$, i.e. $R=0$.
\end{proof}

\begin{lem}\label{l:semi}
Let $\rho'=(\rho,P)$ be the representation of $\wdk$ associated to the natural $l$-adic representation of $\gk$ on $V_l(G)^*$. Then $P=0$ and
$$
\rho\cong\ka\oplus(\chi\otimes\om^{-1}).
$$
\end{lem}

\begin{proof}
Sequence \eqref{eq:111} induces an exact
sequence of corresponding representations of $\wdk$, i.e.
\bbe\label{eq:112}
0\longrightarrow V_l(A)^*\otimes_{\imath}\bC\rarab{h}
V_l(G)^*\otimes_{\imath}\bC\rarab{g} V_l(T)^*\otimes_{\imath}\bC
\longrightarrow 0
\ee
is an exact sequence of $\wdk$-modules, where
$\imath:\bQ_l\hookrightarrow\bC$ is a field embedding, $(\ka,0)$ is the representation of $\wdk$ on $V_l(A)^*\otimes_{\imath}\bC$, $\rho'=(\rho,P)$ is the representation of $\wdk$ on $V_l(G)^*\otimes_{\imath}\bC$, and by Lemma \ref{l:ator}, $(\chi\otimes\om^{-1},0)$ is the representation of $\wdk$ on $V_l(T)^*\otimes_{\imath}\bC$. In particular, \eqref{eq:112} is an exact sequence of $\wk$-modules and it splits if $\rho$ is semisimple, which implies that $\rho\cong\ka\oplus(\chi\otimes\om^{-1})$. Thus, it is enough to show that $P=0$ and $\rho$ is semisimple.

It is known that there is a finite Galois extension $L\subset\Kb$ of $K$ such that $T\times_K L$ splits and $A\times_K L$ has good reduction. Since $\rho$ is semisimple if and only if its restriction to a subgroup of finite index is semisimple (Lemma \ref{l:find}) and $\res_{\cW'(\Kb/L)}\rho'=(\res^L_K\rho,P)$ (\cite{r2}, p. 130), to prove that $P=0$ and $\rho$ is semisimple we can assume that $T$ splits over $K$ and $A$ has good reduction over $K$. Then it follows from Lemma \ref{l:ator} that $\chi$ is trivial. Also, since the image of $I$ under $\rho$ is finite, by Lemma \ref{l:find} to prove that $\rho$ is semisimple it is enough to prove that $\rho(\Phi)$ is diagonalizable. 

Taking into account that $\chi$ is trivial, from \eqref{eq:112} we obtain that in a suitable basis $\rho(\Phi)$
has the following form:
\bbe\label{eq:re}
{\rho(\Phi)= \left(\begin{array}{cc}
\ka (\Phi) & * \\
0 & qE_r \\
\end{array} \right),}
\ee
where $E_r$ is the $r\times r$-identity matrix. Let
$$
\ka _l:\gk\rar \text{Aut}(T_l(A))
$$
be the $l$-adic representation corresponding to the Galois module
$T_l(A)$. It is known that the absolute values of the eigenvalues
of $\ka _l(\Phi)$ are equal to $q^{1/2}$ (\cite{st}, Corollary on
p. 499). Then the absolute values of the eigenvalues of
$\ka(\Phi)$ are equal to $q^{-1/2}$, since the eigenvalues of $\ka
(\Phi)$ are the inverses of the eigenvalues of $\ka _l(\Phi)$. It
follows that none of the eigenvalues of $\ka (\Phi)$ is equal to
$q$. Since $\ka(\Phi)$ is diagonalizable by Lemma \ref{l:semisimple1}, \eqref{eq:re} shows that $\rho(\Phi)$ is diagonalizable,
hence $\rho$ is semisimple and $\rho'$ is admissible.

Let us show now that $P=0$. Since \eqref{eq:112} is an exact sequence of $\wdk$-modules, we have
\begin{eqnarray*}
Nh(x)&=&0,\quad\forall x\in V_l(A)^*\otimes_{\imath}\bC,\\
g(Ny)&=&0, \quad\forall y\in V_l(G)^*\otimes_{\imath}\bC,
\end{eqnarray*}
which implies that $N^2=0$. On the other hand, since $\rho'$ is admissible, it has the following form:
$$
\rho'\cong\bigoplus_{i=1}^s\al_i\otimes\spin(n_i),
$$
where each $\al_i$ is a representation of $\wk$, each $n_i$ is a positive integer, and we can assume that $n_i\ne n_j$ if $i\ne j$ (\cite{r2}, p. 133, Cor. 2). Since $N^2=0$, it follows that each $n_i$ is 1 or 2, i.e. without loss of generality we can assume that
\bbe\label{eq:rho'}
\rho'\cong\al_1\oplus(\al_2\otimes\spin(2)).
\ee

We will show that $\al_2=0$. Assume that $\al_2\ne 0$. 
From \eqref{eq:rho'} we have
$$
\rho\cong\al_1\oplus\al_2\oplus(\al_2\otimes\om). 
$$
On the other hand, since $\rho$ is semisimple, the exact sequence \eqref{eq:112} of $\wk$-modules splits, i.e.
$$
\rho\cong\ka\oplus(\om^{-1})^{\oplus r},
$$
because $\chi$ is trivial. Thus, combining the last two congruences, we get
\bbe\label{eq:113}
\al_1\oplus\al_2\oplus(\al_2\otimes\om)\cong\ka\oplus(\om^{-1})^{\oplus
r}.
\ee
By assumption, $A$ has good reduction, hence by the criterion of
N\'eron-Ogg-\v Shafarevi\v c (\cite{st}, p. 493, Thm. 1) the
inertia group $I$ acts trivially on $V_l(A)^*$. Since by Lemma
\ref{l:semisimple1}, $\ka$ is semisimple it implies that
$\ka=\bigoplus_{i=1}^{2m}\ka_i$, where $m=\dim A$ and
$\ka_1,\ldots,\ka_{2m}$ are one-dimensional subrepresentations of
$\ka$. Thus, it follows from \eqref{eq:113} that $\al_2$ is a sum of
one-dimensional representations. Let $\al_0$ be one of them. Then
using the uniqueness of decomposition of a semisimple module into
simple modules we have from \eqref{eq:113}:
$$
\al_0\cong\om^{-1}\quad\text{or}\quad \al_0\cong\ka_i
$$
for some $\ka_i$, hence $$\al_0\otimes\om\cong
1\quad\text{or}\quad \al_0\otimes\om\cong\ka_i\otimes\om.$$ In
particular, the absolute value of $\al_0\otimes\om (\Phi)$ is $1$
or $q^{-3/2}$, because the absolute value of $\ka_i(\Phi)$ is
$q^{-1/2}$ for each $i$ (see above). It implies that
$\al_0\otimes\om$ is neither $\om^{-1}$ nor $\ka_i$ for any $i$
which contradicts \eqref{eq:113}. Thus, $\al_2=0$ and $\rho'$ is a
representation of $\wk$.
\end{proof}

\begin{lem}\label{l:admi}
$\sg'$ is admissible.
\end{lem}

\begin{proof}
There is the following exact $\gk$-equivariant sequence (\cite{ray}, p. 312):
\bbe\label{eq:mit}
0\longrightarrow G(\Kb)_{l^n}\longrightarrow M(\Kb)_{l^n}\rar
Y(\Kb)/l^nY(\Kb)\longrightarrow 0,
\ee
where $G(\Kb)_{l^n}$ (resp. $M(\Kb)_{l^n}$) denotes
$\Hom(\bZ/l^n\bZ,G(\Kb))$ (resp. $\Hom(\bZ/l^n\bZ,M(\Kb))$). Since
$Y(\Kb)$ is a free group of rank $r$ and $G(\Kb)$ is divisible by
Lemma \ref{l:ediv} (see Appendix \ref{ap:c}), we have the following exact sequence of $\gk$-modules:
$$
0\longrightarrow T_l(G)\longrightarrow T_l(M)\rar
\chi\otimes\bZ_l^r \longrightarrow 0.
$$
By tensoring the above sequence with $\bQ_l$ over $\bZ_l$ and
taking duals over $\bQ_l$ afterwards, we get:
\bbe\label{eq:114}
0\longrightarrow \chi\otimes\bQ_l^r\longrightarrow V_l(M)^*\rar
V_l(G)^* \longrightarrow 0,
\ee
because $\chi\cong\chi^*$ as a representation with finite image,
realizable over $\bZ$. 

As in the proof of Lemma \ref{l:semi}, by Lemma \ref{l:find}
we can assume that $A$ has good reduction over $K$ and $T$ splits over $K$. Then it follows from Lemma \ref{l:ator} that
$\chi$ is trivial.
Also, by Lemma \ref{l:find} to prove that $\sg'$ is admissible it is enough to prove that $\sg(\Phi)$ is diagonalizable.

Sequence \eqref{eq:114} induces an exact sequence of corresponding representations of $\wdk$, i.e.
\bbe\label{eq:114a}
0\longrightarrow (\chi\otimes\bQ_l^r)\otimes_{\imath}\bC\longrightarrow
V_l(M)^*\otimes_{\imath}\bC\rar V_l(G)^*\otimes_{\imath}\bC
\longrightarrow 0
\ee
is an exact sequence of $\wdk$-modules, where $\chi$ is the representation of $\wdk$ on $(\chi\otimes\bQ_l^r)\otimes_{\imath}\bC$, $\sg'=(\sg,N)$ is the representation of $\wdk$ on $V_l(M)^*\otimes_{\imath}\bC$, and by Lemma \ref{l:semi}, $\ka\oplus(\chi\otimes\om^{-1})$ is the representation of $\wdk$ on $V_l(G)^*\otimes_{\imath}\bC$.
Taking into account that $\chi$ is trivial and \eqref{eq:114a} is an exact sequence of $\wk$-modules, we obtain that in a suitable basis $\sg(\Phi)$ has the following form:
\bbe\label{eq:klt}
{\sg(\Phi)= \left(\begin{array}{ccc}
E_r & * & *\\
0 & qE_r & * \\
0 & 0 & \ka (\Phi)
\end{array} \right).}
\ee
Here $\ka(\Phi)$ is diagonalizable by Lemma \ref{l:semisimple1}.
Since the absolute values of the eigenvalues of $\ka(\Phi)$ are equal to $q^{-1/2}$ (see above), none of the eigenvalues of $\ka(\Phi)$ is equal to $1$ or $q$. Thus, \eqref{eq:klt} shows that $\sg(\Phi)$ is diagonalizable, hence $\sg$ is semisimple, and
$\sg'$ is admissible. 
\end{proof}

\begin{proof}[Proof of Proposition \ref{pr:exact_sequence}] 
Since $\sg'$ is admissible by Lemma \ref{l:admi} and the representations of the Weil-Deligne group $\wdk$ on $(\chi\otimes\bQ_l^r)\otimes_{\imath}\bC$ and $V_l(G)^*\otimes_{\imath}\bC$ are actually representations of the Weil group $\wk$, the same argument as in the proof of Lemma \ref{l:semi} applied to \eqref{eq:114a} gives that $\sg'$ has the following form:
\bbe\label{eq:lal}
\sg'\cong\ga\oplus(\de\otimes\spin(2)),
\ee
where $\ga$ and $\de$ are representations of $\wk$. Hence
$$
\sg\cong\ga\oplus\de\oplus(\de\otimes\om).
$$
On the other hand, since $\sg$ is semisimple by Lemma \ref{l:admi}, the exact sequence \eqref{eq:114a} of $\wk$-modules splits, i.e.
$$
\sg\cong\chi\oplus\ka\oplus(\chi\otimes\om^{-1}).
$$
Thus, combining the last two congruences, we get
\bbe\label{eq:117}
\ga\oplus\de\oplus(\de\otimes\om)\cong\chi\oplus\ka\oplus(\chi\otimes\om^{-1}).
\ee
Note that $\chi$ is isomorphic to a subrepresentation of $\ga\oplus(\delta\otimes\om)$, because by \eqref{eq:114a}
$$
\chi\hookrightarrow\ker N=\ga\oplus(\delta\otimes\om).
$$
Thus, $\de$ is isomorphic to a subrepresentation of $\ka\oplus(\chi\otimes\om^{-1})$ by the uniqueness of decomposition of a semisimple module into simple modules. We claim that $\de$ is isomorphic to a
subrepresentation of $\chi\otimes\om^{-1}$. Indeed, suppose there is an irreducible subrepresentation $\de_0$ of $\de$ which is isomorphic to a subrepresentation of $\ka$. Since the absolute values of the eigenvalues of $\ka(\Phi)$ are equal to $q^{-1/2}$ (see above), the eigenvalues of $\de_0(\Phi)$ are of absolute value $q^{-1/2}$. Hence the eigenvalues of $\de_0\otimes\om(\Phi)$ are of absolute value $q^{-3/2}$. On the other hand, it follows from \eqref{eq:117} that $\de_0\otimes\om$ is isomorphic to a subrepresentation of $\chi$, $\ka$, or $\chi\otimes\om^{-1}$, which is a contradiction because the eigenvalues of 
$\chi(\Phi)$, $\ka(\Phi)$, and $\chi\otimes\om^{-1}(\Phi)$ are of absolute values $1$, $q^{-1/2}$, and $q$ respectively. Thus, $\de$ is isomorphic to a
subrepresentation of $\chi\otimes\om^{-1}$. Since $\dim \de=r$ by Lemma
\ref{l:dimens} (see Appendix \ref{ap:c}), we have $\de\cong\chi\otimes\om^{-1}$, hence $\ga\cong\ka$ by \eqref{eq:117} and $\sg'\cong\ka\oplus
(\chi\otimes\om^{-1}\otimes\spin(2))$ by \eqref{eq:lal}.
\end{proof}

\begin{cor}\label{c:kasympl}
$\ka\otimes\om^{1/2}$ is symplectic.
\end{cor}

\begin{proof}
It follows from Proposition \ref{pr:exact_sequence} that
$$
\sg\cong\ka\oplus(\chi\otimes\om^{-1})\oplus\chi.
$$
Thus,
$$
\sg\otimes\om^{1/2}\cong(\ka\otimes\om^{1/2})\oplus
(\chi\otimes\om^{-1/2})\oplus(\chi\otimes\om^{1/2}).
$$
Here $\chi$ is a representation of finite image realizable over
$\bZ$, hence $(\chi\otimes\om^{1/2})^*\cong\chi\otimes\om^{-1/2}$
and $(\chi\otimes\om^{-1/2})\oplus(\chi\otimes\om^{1/2})$ is
symplectic. Since $\sg\otimes\om^{1/2}$ is symplectic by
assumption (see Subsection \ref{sub:loc gen}), this implies that $\ka\otimes\om^{1/2}$ is symplectic
too. Indeed, since $\sg\otimes\om^{1/2}$ is semisimple and symplectic, it follows from Lemma \ref{l:sympl} (see Appendix \ref{sec:b}) that $\sg\otimes\om^{1/2}$ has the following form:
\bbe\label{eq:w1}
\sg\otimes\om^{1/2}\cong\nu\oplus\nu^*\oplus\la_1\oplus\dotsb\oplus\la_s,
\ee
where $\nu$ is a representation of $\wk$ and $\la_1,\ldots,\la_s$ are irreducible symplectic representations of $\wk$. On the other hand, we have
\bbe\label{eq:w2}
\sg\otimes\om^{1/2}\cong(\ka\otimes\om^{1/2})\oplus
(\chi\otimes\om^{1/2})^*\oplus(\chi\otimes\om^{1/2}).
\ee
Putting \eqref{eq:w1} and \eqref{eq:w2} together and taking into account the uniqueness of decomposition of a semisimple module into simple modules, we conclude that $\ka\otimes\om^{1/2}$ must have the same form as $\sg\otimes\om^{1/2}$ in \eqref{eq:w1}, i.e. without loss of generality we can assume that
$$
\ka\otimes\om^{1/2}\cong\mu\oplus\mu^*\oplus\la_1\oplus\dotsb\oplus\la_t,
$$
where $\mu$ is a subrepresentation of $\nu$ and $t\leq s$. Thus, $\ka\otimes\om^{1/2}$ is symplectic.
\end{proof} 

\begin{cor}
The representation $\sg'$ does not depend on choice of $l$ and $\imath$.
\end{cor}
\begin{proof}
The statement is a consequence of Corollary \ref{c:dep} and Proposition \ref{pr:exact_sequence}.
\end{proof}

\begin{cor}\label{cor:obs}
Let $\tau$ be a representation of $\gk$ with real-valued character. Then
\bbe\label{eq:obs1}
W(\sg'\otimes\tau)=W(\ka\otimes\tau)\cdot
\det\tau(-1)^{r}\cdot\det\chi(-1)^{\dim\tau}\cdot
(-1)^{\langle\chi,\tau\rangle}.
\ee
Moreover, when $p>2g+1$ we have
\bbe\label{eq:obs2}
W(\sg'\otimes\tau)=\det\mu(-1)^{\dim\tau}\cdot\det\chi(-1)^{\dim\tau}\cdot
\det\tau(-1)^{r+l_1}\cdot\al^{\dim\tau}\cdot
(-1)^{\langle\chi,\tau\rangle+l_2},
\ee
where $l_1=\dim\mu+\frac{1}{2}(\dim\mu_1+\dotsb+\dim\mu_a)$,
$\mu$ is a representation of $\wk$,
$\mu_1,\ldots,\mu_a$ are irreducible symplectic subrepresentations
of $\ka\otimes\om^{1/2}$ with finite images, $\al=\pm 1$, $l_2=a\cdot\langle 1,\tau\rangle+
a\cdot\langle\eta,\tau\rangle+a\cdot\langle
\hat\mu_1\oplus\dotsb\oplus\hat\mu_a,\tau\rangle$, $\hat{\mu}_1,\ldots,\hat{\mu}_a$ are representations with finite images given by \eqref{eq:hat} such
that $\hat{\mu}_1\oplus\dotsb\oplus\hat{\mu}_a$ is realizable over
$\bQ$, and $\eta$ is
the unramified quadratic character of $K^{\times}$.
\end{cor}
\begin{proof}
Since the root number of a direct sum of
representations of $\wdk$ equals the product of the root numbers
of the summands, we get from Proposition \ref{pr:exact_sequence}
\bbe\label{eq:red}
W(\sg'\otimes\tau)=W(\ka\otimes\tau)\cdot
W(\chi\otimes\om^{-1}\otimes\spin(2)\otimes\tau),
\ee
where by Proposition 6 (\cite{r1}, p. 327)
$$
W(\chi\otimes\om^{-1}\otimes\spin(2)\otimes\tau)=\det\tau(-1)^{r}\cdot\det\chi
(-1)^{\dim\tau}\cdot
(-1)^{\langle\chi,\tau\rangle},
$$
which proves \eqref{eq:obs1}. 

Formula \eqref{eq:obs2} is a consequence of \eqref{eq:obs1} together with Proposition \ref{pr:gen of Rohrlich} and Corollary \ref{c:kasympl}.
\end{proof}

\section{Proof of Theorem \ref{th:main}}\label{sec:th}
We keep the notation of the introduction.
\begin{lem}\label{l:archimed}
Let $X$ be a smooth projective curve of genus $g$ over a number
field $F$ and $\tau$ a representation of $\Gal(\overline F/F)$
with real-valued character. Then at every infinite place $v$ of
$F$ we have
$$
W(X_v,\tau_v)=(-1)^{g\dim\tau}.
$$
\end{lem}
\begin{proof}
To define $W(X_v,\tau_v)$ let $\sg'_v$ denote the representation
of the Weil-Deligne group $\wdl$ associated to $X_v$, then
$W(X_v,\tau_v)=W(\sg'_v\otimes\tau_v)$, where $\tau_v$ is viewed
as a representation of $\wdl$. If $v$ is an infinite place such
that $F_v\cong\bC$, then the representation $\sg'_v=\sg_v$ of
$\cW'(\bC/\bC)=\cW(\bC/\bC)=\bC^{\times}$ has the following form:
$$
\sg_v=(\varphi_{1,0}\otimes H^{1,0})\oplus(\varphi_{0,1}\otimes
H^{0,1}),
$$
where $\varphi_{p,q}:\cW(\bC/\bC)\rar\bC^{\times}$ $(p,q\in\bZ)$
are given by
$$
\varphi_{p,q}(z)=z^{-p}\bar{z}^{-q},
$$
$H^{1,0}$ and $H^{0,1}$ are the components of $H^1(X_v,\bC)$ in
the Hodge decomposition
$$
H^1(X_v,\bC)=H^{1,0}\oplus H^{0,1}.
$$
Here $H^{1,0}$ and $H^{0,1}$ are endowed with the trivial  action
of $\cW(\bC/\bC)$, hence
\bbe\label{eq:arch1}
\sg_v=(\varphi_{1,0}\oplus\varphi_{0,1})^{\oplus g}.
\ee
Let $v$ be an infinite place such that $F_v\cong\bR$. We have
$$
\cW'(\bC/\bR)=\cW(\bC/\bR)=\bC^{\times}\cup J\bC^{\times},
$$
where $J^2=-1$ and $JzJ^{-1}=\bar{z}$ for $z\in\bC^{\times}$. Here
$\cW(\bC/\bC)$ is identified with the subgroup $\bC^{\times}$ of
$\cW(\bC/\bR)$. In this case the representation $\sg'_v=\sg_v$ of
$\cW(\bC/\bR)$ associated to $X_v$ has the following form:
$$
\sg_v=\ind^{\bC}_{\bR}\,\varphi_{0,1}\otimes H^{0,1},
$$
where $\ind^{\bC}_{\bR}\,\varphi_{0,1}$ denotes the representation
of $\cW(\bC/\bR)$ induced from $\varphi_{0,1}$. As in the complex
case, $H^{0,1}$ is endowed with the trivial action of
$\cW(\bC/\bR)$, hence
\bbe\label{eq:arc2}
\sg_v=(\ind^{\bC}_{\bR}\,\varphi_{0,1})^{\oplus g}
\ee
(\cite{r2}, p. 155, \S 20).

It follows from the proof of Theorem 2(i) (\cite{r1}, p. 329) that
\begin{eqnarray*}
W((\varphi_{1,0}\oplus\varphi_{0,1})\otimes\tau_v)&=&(-1)^{\dim\tau}\quad
\text{if}\quad F_v\cong\bC\quad\text{and}  \\
W((\ind^{\bC}_{\bR}\,\varphi_{0,1})\otimes\tau_v)&=&(-1)^{\dim\tau}\quad
\text{if}\quad F_v\cong\bR.
\end{eqnarray*}
Now the statement follows from these formulas together with
formulas \eqref{eq:arch1} and \eqref{eq:arc2}.
\end{proof}

\begin{lem}[\cite{r1}, Lemma on p. 347]\label{l:tau sympl}
Let $G$ be a finite group, $D\subseteq G$ an abelian subgroup, and
$\tau$ an irreducible representation of $G$ with real-valued
character. If $m_{\bQ}(\tau)=2$ then $\res^G_D\tau$ is symplectic.
\end{lem}

\begin{lem}\label{l:dim and det}
Let $G$ be a finite group and $\tau$ an irreducible representation
of $G$ with real-valued character. If $m_{\bQ}(\tau)=2$ 
then $\dim\tau$
is even and $\det\tau$ is trivial.
\end{lem}
\begin{proof}
By Lemma on p. 339 in \cite{r1} if $\tau$ has odd dimension or
nontrivial determinant, then there is a cyclic subgroup $D$ of $G$
such that $\res^G_D\tau$ is not symplectic, which contradicts
Lemma \ref{l:tau sympl}.
\end{proof}
\begin{proof}[Proof of Theorem \ref{th:main}]
By Lemmas \ref{l:archimed} and \ref{l:dim and det},
$W(X_v,\tau_v)=1$ at every infinite place $v$ of $F$.

Let $v$ be a finite place of $F$ lying over a prime number $p$.
Let $\sg_v'$ be the representation of $\wdl$ associated to $X_v$.
Since by Lemma
\ref{l:dim and det} $\det\tau$ is trivial and $\dim\tau$ is even, \eqref{eq:obs1} implies
\bbe\label{eq:11}
W(X_v,\tau_v)=W(\ka_v\otimes\tau_v)\cdot(-1)^{\langle\chi_v,\tau_v\rangle},
\ee
where $\chi_v$ is a representation of $\Gal(\Fbv/F_v)$ realizable
over $\bZ$ (see Subsection \ref{sec:gen} for the definition of
$\chi_v$). Moreover, when $p>2g+1$ from \eqref{eq:obs2} we have
\bbe\label{eq:22}
W(X_v,\tau_v)=(-1)^{a\langle 1,\tau_v\rangle+
a\langle\eta_v,\tau_v\rangle+a\langle{\la},\tau_v\rangle+\langle\chi_v,\tau_v\rangle},
\ee
where $\eta_v$ is the unramified quadratic  character of $F_v^{\times}$, and $\la=\hat\mu_1\oplus\dotsb\oplus\hat\mu_a$ is a representation of
$\Gal(\Fbv/F_v)$ realizable over $\bQ$.

The rest of the proof is analogous to the argument given by D.
Rohrlich in \cite{r1}. Let $K\subset\overline{F}$ be a finite
Galois extension of $F$ such that $\tau$ factors through the group
$G=\Gal(K/F)$ and $\chi_v$ factors through the decomposition
subgroup $H$ of $G$ at $v$. Then
$$
\langle\chi_v,\tau_v\rangle=\langle\ind^G_H\chi_v,\tau\rangle
$$
by Frobenius reciprocity. Since $\chi_v$ is realizable over
$\bQ$, $\ind^G_H\chi_v$ is realizable over $\bQ$, hence
$\langle\ind^G_H\chi_v,\tau\rangle$ is divisible by $m_{\bQ}(\tau)$.
By assumption $m_{\bQ}(\tau)=2$, hence
$\langle\chi_v,\tau_v\rangle$ is even.
Analogously, $\langle 1,\tau_v\rangle$, $\langle\eta_v,\tau_v\rangle$, and $\langle{\la},\tau_v\rangle$ are even, hence $W(X_v,\tau_v)=W(\ka_v\otimes\tau_v)$ by
\eqref{eq:11} and when $p>2g+1$ we have $W(X_v,\tau_v)=1$ by \eqref{eq:22}. When $p\leq 2g+1$ by assumption the decomposition subgroup of $\Gal(L/F)$ at $v$ is abelian, hence $\tau_v$ is symplectic by Lemma \ref{l:tau sympl}. Also
$\ka_v\tens{\om_v}^{1/2}$ is symplectic by Corollary
\ref{c:kasympl}. Since $\ka_v$ is a representation of $\wel$ (see
Subsection \ref{sus: potential}),
$$
W(\ka_v\otimes\tau_v)=W(\ka_v\otimes{\om_v}^{1/2}\otimes\tau_v)=1
$$
by Proposition 2 and the remark after it on p. 319 in \cite{r1}.
\end{proof}

\appendix

\section{}\label{sec:b}
\begin{lem}\label{lem:ir of semid prod}
Let $C=(c)$ be an infinite cyclic group generated by an element
$c$ and let $B=(b)$ be a finite cyclic group of order $n$ generated by an
element $b$. Let  $G=B\rtimes C$ be a semi-direct product, where
$C$ acts on $B$ via $c^{-1}bc=b^k$ for some
$k\in(\bZ/n\bZ)^{\times}$. Denote by $s$ the order of $k$ in
$(\bZ/n\bZ)^{\times}$. Then
every irreducible representation $\la$ of $G$ has the following
form:
$$
\la=\la_0\otimes\phi,
$$
where $\la_0$ is an irreducible representation of $G$ trivial on
the subgroup of $C$ generated by $c^s$ and $\phi$ is a
one-dimensional representation of $G$.
\end{lem}
\begin{proof}
Since $c^s$ is contained in the center of $G$ and $\la$ is an
irreducible complex representation, by Schur's lemma $\la(c^s)$ is
a scalar multiple of the identity matrix $E$, i.e.
$\la(c^s)=a\cdot E$ for some $a\in\bC^{\times}$. Define a
one-dimensional representation $\phi$ of $G$ as follows:
$\phi(b)=1$ and $\phi(c)$ equals an $s$-th root of $a$. Then
$\la_0=\la\otimes\phi^{-1}$ is trivial on $(c^s)$ and
$\la=\la_0\otimes\phi$.
\end{proof}

\begin{proof}[Proof of Proposition \ref{p:sympl irred}]
Let $\la$ be an irreducible symplectic representation
of $G$. Then by Lemma \ref{lem:ir of semid prod},
$\la=\la_0\otimes\phi$, where $\la_0$ is an irreducible
representation of $G$ trivial on the subgroup of $C$ generated by
$c^s$ and $\phi$ is a one-dimensional representation of $G$. Since
$\la$ is symplectic, $\la$ and its contragredient representation
have the same character, which implies that for any $g\in G$ we
have
$$
\phi(g)\cdot\text{tr}\,\la_0(g)=\phi(g)^{-1}\cdot\text{tr}\,\la_0(g^{-1}).
$$
\noindent
Taking into account that $\la_0$ is trivial on $(c^s)$, the above
equation for $g=c^s$ gives $\phi(c^{2s})=1$, i.e. $\la$ can be
considered as an irreducible symplectic representation of the
finite group $H=G/(c^{2s})\cong B\rtimes C/(c^{2s})$. By abuse of
notation we will denote the image of $c$ in $C/(c^{2s})$ also by
$c$, then $c^{2s}=1$ and $c^{-1}bc=b^k$ in $H$. As an irreducible
representation of the semi-direct product $H$, $\la$ can be constructed from a one-dimensional representation $\psi_1$ of $B$ in the following way. Let
$\psi_1(b)=\xi$ for some $n$-th root of unity $\xi$ of order $d$
in $\bC^{\times}$. Let $\Ga=(c^x)$, where $x=|k|$ in
$(\bZ/d\bZ)^{\times}$, and $\psi_2$ be a one-dimensional
representation of $\Ga$. Then $\psi_1$ and $\psi_2$ can be extended
to representations of $B\rtimes\Ga$ via
\begin{eqnarray*}
\psi_1(c^{xv}b^t)&=&\psi_1(b^t), \\
\psi_2(c^{xv}b^t)&=&\psi_2(c^{xv}).
\end{eqnarray*}
Then $\la=\text{Ind}_{B\rtimes\Ga}^H(\psi_1\otimes\psi_2)$
(\cite{S}, p. 62, Prop. 25). Let $W$ be a representation space of
$H$ corresponding to $\la$, $V=\bC e \subseteq W$ be a
one-dimensional subrepresentation of $\text{Res}^H_{B\rtimes
\Ga}\la$ isomorphic to $\psi_1\otimes\psi_2$ and spanned by a
nonzero vector $e\in V$ over $\bC$. Then $W=V\oplus cV\oplus
c^2V\oplus\cdots\oplus c^{x-1}V$ and $\la$ has the following form
in the basis $\{e,ce,c^{2}e,\ldots,c^{x-1}e\}$:
\begin{displaymath}
{\la(b)}= \left(\begin{array}{ccccc}
\xi & 0 & 0 & \ldots & 0 \\
0 & \xi^k & 0 & \ldots & 0 \\
0 & 0 & \xi^{k^2} & \ldots & 0 \\
\vdots & \vdots & \vdots & \ddots & \vdots \\
0 & 0 & 0 & \ldots & \xi^{k^{x-1}} \\
\end{array} \right), \quad
{\la(c)}= \left(\begin{array}{ccccc}
0 & 0 & 0 & \ldots & \psi_2(c^x) \\
1 & 0 & 0 &\ldots & 0 \\
0 & 1 & 0&\ldots & 0 \\
\vdots & \vdots & \vdots & \ddots & \vdots \\
0 & 0 & 0 & \ldots & 0 \\
\end{array} \right).
\end{displaymath}
\noindent
Since $\la$ is symplectic, $x=\dim\la$ is even and $\det\la=1$,
hence $\det\la(c)=-\psi_2(c^x)=1$, which implies that
$\psi_2(c^x)=-1$. Denote by $\chi$ the character of $\la$. By
Proposition 39 (\cite{S}, p. 109), $\la$ is symplectic if and only
if
\begin{equation}\label{eq:sym}
-1=\frac{1}{|H|}\cdot\sum_{y\in H}\chi(y^2).
\end{equation}
\noindent
Let $y=c^vb^t$, consequently, $y^2=c^{2v}b^{t(1+k^v)}$. Clearly,
$\chi(y^2)=0$ if $y^2\not\in B\rtimes\Ga$ and $y^2\in B\rtimes\Ga$
if and only if $x$ divides $2v$ and, since $x$ is even, if and
only if $\frac{x}{2}$ divides $v$. Let $v=\frac{x}{2}m$, then we
have
\begin{equation}
\begin{split}
\sum_{y\in H}\chi(y^2)=\sum_{ \substack{
y\in H\\
y^2\in B\rtimes\Ga }}
\chi(y^2)=\sum_{t,m}\chi(c^{mx}b^{t(1+k^v)}) \\
=\sum_{\substack{m,t \\ m\text{ even}}}\chi(b^{t(1+k^v)})-
\sum_{\substack{m,t \\ m\text{ odd}}}\chi(b^{t(1+k^v)}).
\end{split}
\end{equation}
Let $S_1=\sum\limits_{m \text{
even}}\sum\limits_{t}\chi(b^{t(1+k^v)})$ and $S_2=\sum\limits_{m
\text{ odd}}\sum\limits_{t}\chi(b^{t(1+k^v)})$.

If $m$ is even, then $v=x(\frac{m}{2})$ and, since $x=|k|$ in 
$(\bZ/d\bZ)^{\times}$, $1+k^v\equiv 2\,(\text{mod}\,d)$. Since $\chi
(b^t)=\xi^t+\xi^{kt}+\cdots +\xi^{k^{x-1}t}$ and $\xi^d=1$, we
have $\chi(b^{t(1+k^v)})=\chi(b^{2t})$ and
$S_1=\sum\limits_{m\text{ even}}\sum\limits_{t}\chi(b^{2t})$. We
will show that $\sum\limits_{t}\chi (b^{2t})=0$. First, note that
if $d=1,2$, then $\la$ is one-dimensional, hence cannot be
symplectic. 

If $r\in\bZ$ and $r\equiv0\,(\text{mod}\,d)$ then
$$
\sum_{t=0}^{n-1}\chi(b^{rt})=\sum_{t=0}^{n-1}\sum_{j=0}^{x-1}\xi^{rtk^j}=nx.
$$

If $r\in \bZ$ and $r\not\equiv0\,(\text{mod}\,d)$ then
$$
\sum_{t=0}^{n-1}\chi(b^{rt})=\sum_{t=0}^{n-1}\sum_{j=0}^{x-1}\xi^{rtk^j}=
\sum_{j=0}^{x-1}\frac{1-\xi^{rnk^j}}{1-\xi^{rk^j}}=0.
$$
\noindent
Thus
\begin{eqnarray}\label{e:chi}
{\sum_{t=0}^{n-1}\chi(b^{rt})}= \left\{\begin{array}{cc}
nx, & r\equiv0\,(\text{mod}\,d); \\
0, & r\not\equiv0\,(\text{mod}\,d).
\end{array}\right.
\end{eqnarray}
\noindent
Since $d\not= 1,2$, formula \eqref{e:chi} implies that $S_1=0$.

If $m$ is odd, then $k^v\equiv k^{\frac{x}{2}}\,(\text{mod}\,d)$,
hence $\chi(b^{t(1+k^v)})=\chi(b^{t(1+k^{\frac{x}{2}})})$. Thus
$$
S_2=\sum\limits_{m
\text{ odd}}\sum_{t=0}^{n-1}\chi(b^{t(1+k^v)})=\frac{2s}{x}\cdot
\sum_{t=0}^{n-1}\chi(b^{t(1+k^{\frac{x}{2}})})
$$
and by \eqref{e:chi} we have
\bbe\label{eq:jk}
{S_2}= \left\{\begin{array}{cc}
2sn, & 1+k^{\frac{x}{2}}\equiv0\,(\text{mod}\,d); \\
0, & 1+k^{\frac{x}{2}}\not\equiv0\,(\text{mod}\,d).
\end{array}\right.
\ee
\noindent
Hence
$$
\frac{1}{|H|}\cdot\sum_{y\in H}\chi(y^2)=\frac{1}{2sn}\cdot(S_1-S_2)=-\frac{S_2}{2sn}
$$
which together with \eqref{eq:sym} and \eqref{eq:jk}
proves the proposition.
\end{proof}

Let $D$ be a group, $U$ a finite-dimensional $\bC[D]$-module, and
$U^*$ the dual $\bC[D]$-module of $U$. Let $\check U$ denote the
vector space over $\bC$ with the underlying abelian group $U^*$
and multiplication by constants defined as follows:
$$
a\cdot\phi=\overline{a}\phi,\quad a\in\bC,\,\phi\in U^*,
$$
where $\overline{a}$ is the complex conjugate of $a$. Clearly, the
$\bC[D]$-module structure on $U^*$ makes $\check U$ into a
$\bC[D]$-module. In what follows by $\check U$ we mean a
$\bC[D]$-module with this structure. We say that $U$ is {\it
unitary} if $U$ admits a nondegenerate invariant hermitian form
(not necessarily positive definite).

\begin{lem}\label{l:sympl}
Every semisimple unitary, orthogonal, or symplectic representation
$\la$ of a group $D$ has the following form
$$
\la\cong\nu\oplus\tilde\nu\oplus\la_1^{z_1}\oplus\cdots\oplus\la_t^{z_t},
$$
where $\nu$ is a representation of $D$, $\tilde\nu=\nu^*$ if $\la$
is orthogonal or symplectic and $\tilde\nu=\check\nu$ if $\la$ is
unitary, $\la_1,\ldots,\la_t$ are pairwise nonisomorphic
irreducible unitary, orthogonal, or symplectic representations of
$D$ respectively.
\end{lem}

\begin{proof}[Proof of Lemma \ref{l:sympl}]
We say that a unitary, orthogonal, or symplectic representation is
{\it minimal} if it cannot be written as an orthogonal sum of
nonzero invariant subspaces. Clearly, every unitary, orthogonal,
or symplectic representation is an orthogonal sum of minimal unitary,
orthogonal, or symplectic representations respectively. Thus, it
is enough to prove that if $\la$ is a semisimple minimal unitary,
orthogonal, or symplectic representation of $D$, then either $\la$
is irreducible or $\la\cong\nu\oplus\tilde\nu$ for some
irreducible representation $\nu$ of $D$. Let $U$ be a
representation space of $D$ corresponding to $\la$, $\tilde U=U^*$
if $\la$ is orthogonal or symplectic and $\tilde U=\check U$ if
$\la$ is unitary. Since $\la$ is semisimple,
$U=V_1\oplus\dotsb\oplus V_n$, where $V_1,\ldots ,V_n$ are nonzero
simple $\bC [D]$-submodules of $U$. Let $\langle\cdot\,
,\cdot\rangle$ be a nondegenerate invariant form on $U$. It
defines a $\bC[D]$-module isomorphism $\phi$ between $U$ and
$\tilde U$ via $\phi(u)=\langle u\, ,\cdot\rangle$, $u\in U$. Let
$\psi:\tilde U\longrightarrow
\tilde{V_1}\oplus\dotsb\oplus\tilde{V_n}$ denote the usual
isomorphism between $\tilde U={(V_1\oplus\dotsb\oplus V_n)}^{\sim}$ and
$\tilde{V_1}\oplus\dotsb\oplus\tilde{V_n}$. For each $i$ and $j$
let $\al_{ij}:V_i \longrightarrow\tilde{V_j}$ be a $\bC[D]$-module
homomorphism defined by the following diagram:
$$
\xymatrix{
U \ar[r]^-{\psi\circ\phi}   &   \tilde{V_1}\oplus\dotsb\oplus \tilde{V_n} \ar[d]^{\pi_j} \\
V_i \ar@{^{(}->}[u] \ar[r]^{\al_{ij}} &  \tilde{V_j} }
$$
where $\pi_j$ is the projection onto $j$-th factor. Since
$\psi\circ\phi$ is an isomorphism, there exists some $\tilde{V_i}$
such that $\al_{1i}\ne 0$, which implies that $\al_{1i}$ is an
isomorphism, since $V_1,\ldots,V_n$ are simple. If $i=1$, then it
follows that $\langle\cdot\, ,\cdot\rangle\vert_{V_1}$ is
nondegenerate, hence $V_1$ and its orthogonal complement are
invariant subspaces of $U$. Since $U$ is minimal, it implies that
$U=V_1$ and $U$ is irreducible. Thus, we can assume that for each
$j$ we have $\al_{jj}=0$, which is equivalent to $\langle V_j
\,,V_j\rangle=0$. Without loss of generality we can assume that
$\al_{12}\ne 0$. Then $\al_{21}\ne 0$. Indeed, if $\al_{12}\ne 0$,
then there is some $u\in V_1$ such that $\langle u
\,,\cdot\rangle\vert_{V_2}\ne 0$, i.e. there is some $v\in V_2$
such that $\langle u \,,v\rangle\ne 0$, hence $\langle v
\,,u\rangle\ne 0$, which is equivalent to $\al_{21}\ne 0$. Let us
prove now that $\langle\cdot\, ,\cdot\rangle\vert_{V_1\oplus V_2}$
is nondegenerate. Let $u+v\in V_1\oplus V_2$ and $\langle u+v \,
,x+y\rangle=0$ for any $x+y\in V_1\oplus V_2$. We have $\langle
u+v \, ,x+y\rangle=\langle u\, ,y\rangle +\langle v \,
,x\rangle=0$, because $\langle V_1 \, ,V_1\rangle=\langle V_2 \,
,V_2\rangle=0$. Take $x=0$ in this equation, then $\langle u \,
,y\rangle=0$ for any $y\in V_2$, hence $u=0$, because
$\al_{12}(V_1)=\tilde{V_2}$. Analogously, $v=0$. Since $U$ is
minimal, the same argument as above implies that $U=V_1\oplus
V_2\cong V_1\oplus\tilde{V_1}$.
\end{proof}


\begin{proof}[Proof of Proposition \ref{p:realiz q}]
By Lemma \ref{l:sympl}
\bbe\label{eq:la}
\la\cong\nu\oplus\nu^*\oplus\la_1^{z_1}\oplus\dotsb\oplus\la_t^{z_t},
\ee
where $\nu$ is a representation of $G$ and $\la_1,\ldots,\la_t$
are pairwise nonisomorphic irreducible symplectic representations
of $G$. Let $\nu=\nu_1^{l_1}\oplus\dotsb\oplus\nu_r^{l_r}$, where
$\nu_1,\ldots,\nu_r$ are pairwise nonisomorphic irreducible representations of $G$. By Lemma
\ref{lem:ir of semid prod} for each $i$ we have
$\nu_i=\nu_i^{0}\otimes\phi_i$, where $\phi_i$ is a
one-dimensional representation of $G$ and $\nu_i^{0}$ is an
irreducible representation of $G$ trivial on $(c^s)$. It follows
that $\nu_i^{0}$ can be considered as a representation of
$H=G/(c^{2s})$ and as a representation of $H$ it can be written in
the following form $\nu_i^{0}=\text{Ind}_{B\rtimes\Ga_i}^H\psi_i$,
where $\psi_i$ is a one-dimensional representation of $B\rtimes\Ga_i$,
$\psi_i(b)=\xi_i$ for an $n$-th root of unity $\xi_i$ of
order $d_i$, $x_i=|k|$ in $(\bZ/d_i\bZ)^{\times}$, and $\Ga_i=(c^{x_i})$
(see the proof of Proposition \ref{p:sympl irred} and \cite{S}, p.
62, Prop. 25). Thus
\begin{displaymath}
{\nu_i^{0}(b)}= \left(\begin{array}{cccc}
\xi_i & 0 & \ldots & 0 \\
0 & \xi_i^{k} & \ldots & 0 \\
\vdots & \vdots & \ddots & \vdots \\
0 & 0 & \ldots & \xi_i^{k^{x_i-1}} \\
\end{array} \right), \end{displaymath}
\begin{displaymath}
{\nu_i(b)=\nu_i^{0}\otimes \phi_i(b)}= \left(\begin{array}{cccc}
\xi_i\phi_i(b) & 0 & \ldots & 0 \\
0 & (\xi_i\phi_i(b))^{k} & \ldots & 0 \\
\vdots & \vdots & \ddots & \vdots \\
0 & 0 & \ldots & (\xi_i\phi_i(b))^{{k}^{x_i-1}} \\
\end{array} \right).
\end{displaymath}
\noindent
In the second matrix we used the relation $\phi_i(b)^{k-1}=1$,
which follows from the fact, that $\phi_i$ is a one-dimensional
representation of $G$ and $c^{-1}bc=b^k$. By Proposition \ref{p:sympl irred} each $\la_i=\text{Ind}^H_{B\rtimes L_i}\rho_i$, where $\rho_i$ is a one-dimensional representation of $B\rtimes L_i$,
$\rho_i(b)=\eta_i$ for an $n$-th root of unity $\eta_i$ of order
$u_i$, $y_i=|k|$ in $(\bZ/u_i\bZ)^{\times}$, $L_i=(c^{y_i})$, and $\rho_i(c^{y_i})=-1$. 

We will need the following lemma:

\begin{lem}\label{l:pol}
Let $d_1,\ldots ,d_m$ be pairwise distinct natural numbers. For
each $d_i$ let $p_i(X)\in\bC[X]$ be a monic polynomial, all the
roots of which are some primitive $d_i$-th roots of unity and let
$p(X)=p_1(X)\cdots p_m(X)$. If $p(X)\in\bQ[X]$, then each $p_i(X)$
is a power of the $d_i$-th cyclotomic polynomial $\Phi_{d_i}(X)$.
\end{lem}
\begin{proof}
Let $\xi_1$ be a root of $p_1(X)$, then $\xi_1$ is a primitive
$d_1$-th root of unity. Since $\xi _1$ is a root of
$p(X)\in\bQ[X]$ and $\Phi_{d_1}(X)$ is the minimal polynomial of
$\xi_1$ over $\bQ$, $p(X)=\Phi_{d_1}(X)^{s_1}\cdot r(X)$, where $s_1$ is
a natural number and $r(X)$ is a polynomial over $\bQ$ not
divisible by $\Phi_{d_1}(X)$. Applying the same argument to all $p_i(X)$
and taking into account that $\Phi_{d_1}(X),\ldots ,\Phi_{d_m}(X)$ are
distinct, we get
$p(X)=\Phi_{d_1}(X)^{s_1}\cdots\Phi_{d_m}(X)^{s_m}$ for some
$s_1,\ldots, s_m>0$. Since the roots of each $p_i(X)$ can be only
primitive $d_i$-th roots of unity, we conclude that
$p_i(X)=\Phi_{d_i}(X)^{s_i}$ for each $i$.
\end{proof}

Since the characteristic polynomial $p$ of $\la(b)$ has
coefficients in $\bQ$, by Lemma \ref{l:pol} we can assume that $\xi_1\phi_1(b),\ldots,\xi_r\phi_r(b),\eta_1,\ldots,\eta_t$
are primitive roots of unity of the same order $d$ and that
$p=\Phi_d^v$ for some $v$, where $\Phi_d$ is the $d$-th
cyclotomic polynomial. Indeed, $\la$ can be written as a sum of semisimple symplectic representations of $G$ which have this property and it is enough to show that for each of them \eqref{eq:realiz1} holds.

Let $x=|k|$ in $(\bZ/d\bZ)^{\times}$ and
$\Ga=(c^x)$. If $\la\cong\nu\oplus\nu^*$ then there is nothing to
prove. Thus, we assume that there is $\la_1$ in \eqref{eq:la}.
Since $\la_1$ is symplectic, $x$ is even, $d\neq 1,2$, and
$k^{\frac{x}{2}}\equiv -1\,(\text{mod}\,d)$ by Proposition
\ref{p:sympl irred}. Note, that $x$ divides each $x_i$. Indeed,
$(\xi_i\phi_i(b))^{k^{x_i}}=\xi_i\phi_i(b)$, hence $k^{x_i}\equiv
1\,(\text{mod}\,d)$, because $\xi_i\phi_i(b)$ is a primitive
$d$-th root of unity by assumption. For each $i$ denote by
$p_{\nu_i}$ the characteristic polynomial of $\nu_i(b)$ and by
$p_{\nu_i^*}$ the characteristic polynomial of $\nu_i^*(b)$. Then
$p_{\nu_i}=p_{\nu_i^*}$. This is true because $x$ divides $x_i$,
$x$ is even, $k^{\frac{x_i}{2}}\not\equiv 1\,(\text{mod}\,d)$, and
$k^{\frac{x}{2}}\equiv -1\,(\text{mod}\,d)$, hence each root of
$p_{\nu_i}$ appears in $p_{\nu_i}$ with its complex conjugate.
Thus
$$
p=p_{\nu_1}^{2l_1}\cdots p_{\nu_r}^{2l_r}p_{\la_1}^{z_1}\cdots
p_{\la_t}^{z_t},
$$
where for each $i$ we denote by $p_{\la_i}$ the characteristic
polynomial of $\la_i(b)$.

For each primitive $d$-th root of unity $\xi$ write
$q(\xi)=(X-\xi)(X-\xi^k)\cdots(X-\xi^{k^{x-1}})$, where $x=|k|$ in $(\bZ/d\bZ)^{\times}$. Clearly, all
$\xi,\xi^k,\ldots,\xi^{k^{x-1}}$ are distinct and for two
primitive $d$-th roots of unity $\xi$ and $\xi'$ either
$q(\xi)=q(\xi')$ or $q(\xi)$ and $q(\xi')$ have no common roots.
In this notation $p_{\la_i}=q(\eta_i)$ and
$p_{\nu_i}=q(\xi_i\phi_i(f))^{\alpha_i}$, where
$\alpha_i=\frac{x_i}{x}$. Since $\la_1,\ldots,\la_t$ are irreducible, symplectic, and pairwise
nonisomorphic, it follows from Proposition \ref{p:sympl irred} that $q(\eta_i)\ne q(\eta_j)$ for $i\ne j$. Without loss
of generality we can assume that $p$ has the following form:
\bbe\label{eq:p}
p=q(\xi_1\phi_1(b))^{2m_1}\dotsb
q(\xi_f\phi_f(b))^{2m_f}q(\eta_1)^{z_1}\dotsb q(\eta_t)^{z_t},
\ee
where $f\leq r$, $m_1,\ldots,m_f$ are positive integers, and
$q(\xi_1\phi_1(b)),\ldots, q(\xi_f\phi_f(b))$ have no common
roots. There are two possibilities:

\begin{enumerate}
\item there exists some $q(\xi_i\phi_i(b))$ which is not equal to any of
$q(\eta_1),\ldots,q(\eta_t)$. Without loss of generality we can
assume that $i=1$;
\item each $q(\xi_i\phi_i(b))$ equals some $q(\eta_j)$.
\end{enumerate}
\bigbreak (1) In this case, since $p=\Phi_d^v$, it follows from
$\eqref{eq:p}$ that for each $j$ we have $z_j+2\cdot\al(j)=2m_1$,
where $\al(j)=m_{\be}$ if $q(\eta_j)$ equals some $q(\xi_{\be}\phi_{\be}(b))$
and $\al(j)=0$ otherwise. Thus, in this case all $z_1,\ldots,z_t$
are even and $[\la]=[\nu]+[\nu^*]+2\cdot[\mu_0]$, where
$\mu_0=\la_1^{\frac{z_1}{2}}\oplus\dotsb\oplus\la_t^{\frac{z_t}{2}}$
is symplectic of finite image because all $\la_1,\ldots,\la_t$ are
symplectic of finite images. 
\bigbreak (2) In this case, since
$p=\Phi_d^v$, it follows from $\eqref{eq:p}$ that for each $j$ we
have $z_j+2\cdot \al(j)=v$, where $\al(j)=m_{\be}$ if $q(\eta_j)$
equals some $q(\xi_{\be}\phi_{\be}(b))$ and $\al(j)=0$ otherwise.
Moreover, it follows that $q(\eta_1)\dotsb q(\eta_t)=\Phi_d$. Thus
$$
[\la]=[\nu]+[\nu^*]-2\cdot[\mu'_0]+v\cdot[\la_1]+\dotsb+v\cdot[\la_t],
$$
where $\mu'_0=\la_1^{\al(1)}\oplus\dotsb\oplus\la_t^{\al(t)}$ is
symplectic of finite image and it is enough to show that
$\hat\la_1\oplus\dotsb\oplus\hat\la_t$ is realizable over $\bQ$.
Recall that for each $i$, $\hat\la_i=\text{Ind}_{B\rtimes
\Ga}^H\varphi_i$, where $\varphi_i(b)=\xi_i$ for some primitive
$d$-th root of unity $\xi_i$, $x=|k|$ in $(\bZ/d\bZ)^{\times}$,
$\Ga=(c^x)$, and $\varphi_i(c^x)=1$ (see Proposition \ref{p:sympl
irred} and \eqref{eq:hat}). Since the representations of this form
are completely defined by a root of unity $\xi$, we will denote
them by $\Theta(\xi)$. For any $r$ dividing $d$ the cyclic group
$(k)$ acts on the set of all primitive $r$-th roots of unity via
$\xi\longmapsto\xi^k$. Let $\{\xi_r^1,\ldots,\xi_r^{w_r}\}$ be the
set of representatives for this action and let
$$
\Theta(r)=\sum_{i=1}^{w_r}\Theta(\xi^i_r).
$$
Then the characteristic polynomial of $\Theta(r)(b)$ is just
$\Phi_r$. Since the characteristic polynomial of
$\hat\la_1\oplus\dotsb\oplus\hat\la_t$ is $q(\eta_1)\dotsb
q(\eta_t)=\Phi_d$, it follows that
$\hat\la_1\oplus\dotsb\oplus\hat\la_t=\Theta(d)$. By induction on
$d$ we will prove that each $\Theta(d)$ is realizable over $\bQ$.

Clearly, $\Theta(d)$ is realizable over $\bQ$ when $d=1$, because
in this case $\Theta(d)=1$. Let $L=(b^d)\rtimes C$ and
$\pi=\text{Ind}_L^H1$. Then the characteristic polynomial of
$\pi(b)$ is $x^d-1$, consequently,
$$
\pi\cong\sum_{r|d}\Theta(r) \quad \text{on} \quad B.
$$
We will prove that this is true on the whole group $H$. Observe,
that for any $r$ all $\Theta(\xi_r^1),\ldots,\Theta(\xi_r^{w_r})$
are irreducible over $\bC$ and
$\Theta(\xi_r^i)\cong\Theta(\xi_{r'}^{i'})$ only if $i=i'$ and
$r=r'$. Let $\chi^i_r$ be the character of $\Theta(\xi^i_r)$.
Then, using Frobenius reciprocity, we have
$$
\langle\pi,\Theta(\xi^i_r) \rangle=\langle\text{Ind}_L^H
1,\Theta(\xi^i_r)\rangle=\langle 1,
\text{Res}^H_L\Theta(\xi^i_r)\rangle=\frac{d}{2ns}\sum_{u,v}\chi^i_r(b^{du}c^v)=
\frac{d}{2ns}\sum_{u,v}\chi^i_r(c^v)=1,
$$
hence $\pi\cong\sum_{r|d}\Theta(r)$ on $H$. Since $\pi$ is
realizable over $\bQ$ and $\Theta(r)$ is realizable over $\bQ$ for
any $r<d$ by induction, $\Theta(d)=\pi-\sum_{r|d,r\ne d}\Theta(r)$
is realizable over $\bQ$.
\end{proof}

\begin{proof}[Proof of Proposition \ref{pr:gen of Rohrlich}]
Let $\la=\ka\tens\om^{1/2}$. Then
$W(\ka\otimes\tau)=W(\la\otimes\tau)$, because real powers of
$\om$ do not change the root number. Since the root number of
representations of $\wk$ is multiplicative in short exact
sequences, there is the unique homomorphism
$$
\al:R(\wk)\longrightarrow\bC^{\times}
$$
such that $\al([\la])=W(\la)$ for any representation $\la$ of
$\wk$. Thus, it follows from Corollary \ref{cor:ad} that
\bbe\label{eq:root}
W(\la\otimes\tau)=W(\mu\otimes\tau)\cdot W(\mu^*\otimes\tau)\cdot
\frac{W(\mu_0\otimes\tau)^2}{W(\mu'_0\otimes\tau)^2}\cdot
W(\mu_1\otimes \tau)\dotsb W(\mu_a\otimes\tau).
\ee
Since $\tau$ has finite image and real-valued character, we have
\begin{eqnarray*}
\notag W(\mu\otimes\tau)\cdot W(\mu^*\otimes\tau)&=&
W(\mu\otimes\tau)\cdot W((\mu\otimes\tau)^*) \\
&=&
\det(\mu\otimes\tau)(-1)=\det\mu(-1)^{\dim\tau}\cdot\det\tau(-1)^{\dim\mu}.
\end{eqnarray*}
Also, since $\mu_0$ and $\mu'_0$ are symplectic and of finite
images, $W(\mu_0\otimes\tau)=\pm 1$, $W(\mu'_0\otimes\tau)=\pm 1$
(\cite{r1}, p. 315), hence from \eqref{eq:root} we get
\bbe\label{eq:w}
W(\la\otimes\tau)=\det\mu(-1)^{\dim\tau}\cdot\det\tau(-1)^{\dim\mu}\cdot
W(\mu_1\otimes \tau)\dotsb W(\mu_a\otimes\tau).
\ee
Thus, we need to compute
$W(\mu_1\otimes\tau),\ldots,W(\mu_a\otimes\tau)$. Let $\gamma$ be
an irreducible symplectic subrepresentation of $\la$. Let
$L/K^{unr}$ be a minimal subextension of $\Kb/K^{unr}$ over which
$A$ acquires good reduction. Then, as was discussed in Subsection
\ref{sus: potential}, $\la$ and, consequently $\ga$, can be
considered as a representation of $G=B\rtimes(\Phi)$, where
$B=\Gal(L/K^{unr})$ is a finite cyclic group, $(\Phi)$ is an
infinite cyclic group. Let $x=\dim\ga$. Then by Proposition
\ref{p:sympl irred}, as a representation of $G$, $\ga$ is induced
from a one-dimensional representation of $B\rtimes(\Phi^{x})$.
Hence, as a representation of $\wk$, $\ga$ is induced from a
one-dimensional representation $\phi$ of $\cW(\Kb/H_{x})$, where
$H_{x}$ is the unramified extension of $K$ of degree $x$, i.e.
$\ga=\ind^{H_{x}}_K\phi$. Since $\ga$ is symplectic, $x$ is even.
Let $x=2y$ and let $H_y$ be the unramified extension of $K$ of degree
$y$, hence $K\subseteq H_y\subseteq H_{x}$. Let
$\ga'=\ind^{H_{x}}_{H_y}\phi$, $\tau'=\res^{H_y}_{K}\tau$, then
$\ga=\ind^{H_y}_K\ga'$ and by Formula (1.4) (\cite{r1}, p. 316) we
have
\bbe\label{eq:10}
W(\ga\otimes\tau)=W(\ind^{H_y}_K(\ga'\otimes\tau'))=W(\ga'\otimes\tau')W(\ind
^{H_y}_K1_{H_y})^{2\dim\tau}.
\ee

Let us prove first that $W(\ind^{H_y}_K1_{H_y})^{2\dim\tau}=1$.
Let $\varpi$ be a uniformizer of $K$. It is easy to check that
$\ind^{H_y}_K1_{H_y}=\bigoplus^{y-1}_{i=0}\chi_i$, where
$\chi_0,\ldots,\chi_{y-1}$ are all the unramified characters of
$K^{\times}$, satisfying $\chi_i(\varpi)^y=1$. Hence
$W(\ind^{H_y}_K1_{H_y})=\prod^{y-1}_{i=0}W(\chi_i)$. By Formula
($\eps 3$) (\cite{r2}, p. 142) for each $i$ we have
$W(\chi_i)=\xi_i^{n(\psi)}$, where $n(\psi)\in\bZ$, each
$\xi_i$ is a $y$-th root of unity, and $\xi_i\ne\xi_j$ if $i\ne j$. Hence
$\prod^{y-1}_{i=0}W(\chi_i)=\prod^{y-1}_{i=0}\xi_i^{n(\psi)}=\pm1$
and
\bbe\label{eq:20}
W(\ind^{H_y}_K1_{H_y})^{2\dim\tau}=1.
\ee

To compute $W(\ga'\otimes\tau')$ we will show that Theorem
\ref{th:rohrlich} can be applied to $H_y$, $\ga'$, and $\tau'$.
Indeed, $\tau'$ is a representation of $\Gal(\Kb/H_y)$ with
real-valued character and $\ga'=\ind^{H_{x}}_{H_y}\phi$ is a
two-dimensional representation of $\cW(\Kb/H_y)$ induced from a
character $\phi$ of finite image (by Proposition \ref{p:sympl
irred}), hence $\ga'$ is a representation of $\Gal(\Kb/H_y)$.
Since $\ind^{H_y}_K\ga'=\ga$ is irreducible, $\ga'$ is irreducible
too. Since $\dim\ga'=2$, $\ga'$ is symplectic if and only if
$\det\ga'$ is trivial, because $\text{Sp}(2,\bC)=\text{SL}(2,\bC)$ (\cite{r1}, p. 317). From Proposition \ref{p:sympl irred} we
find that as a representation of $B\rtimes(\Phi^y)$, $\ga'$ has
the following form
\begin{displaymath}
{\ga'(f)}= \left(\begin{array}{cc}
\xi & 0 \\
0 & \xi^{-1}\\
\end{array} \right), \quad
{\ga'(\Phi^y)}= \left(\begin{array}{cc}
0 & -1\\
1 & 0\\
\end{array} \right),
\end{displaymath}
where $f$ is a generator of $B$, $\xi$ is a root of unity. It
follows immediately that $\det\ga'=1$. Thus to be able to apply
Theorem \ref{th:rohrlich}, we need only to check that $\phi$ is a
tame character of $H_{x}^{\times}$. It follows from the fact that $\phi$ is trivial on $\Gal(\Kb/L)$ and $L/K^{unr}$ is tamely ramified, because $p>2m+1$. We have shown
that $\tau'$ is a representation of $\Gal(\Kb/H_y)$,
$\ga'=\ind^{H_{x}}_{H_y}\phi$ is a two-dimensional symplectic,
irreducible representation of $\Gal(\Kb/H_y)$, and $\phi$ is a
tame character of $H_{x}^{\times}$. By Theorem \ref{th:rohrlich}
\bbe\label{eq:3}
W(\ga'\otimes\tau')=\det\tau'(-1)\cdot\varphi^{\dim\tau}\cdot(-1)^{\langle
1,\tau'\rangle+ \langle\eta',\tau'\rangle+\langle
\hat{\la'},\tau'\rangle},
\ee
where $\eta'$ is the unramified quadratic character of
$H_y^{\times}$,
$\hat{\ga'}=\ind^{H_{x}}_{H_{y}}(\phi\otimes\theta)$, $\theta$ is
the unramified quadratic character of $H_{x}^{\times}$, and
$\varphi=\pm1$.

Since $\tau'=\res^{H_y}_K\tau$, we have $\det\tau'=\det\tau\circ
N_{H_y/K}$, hence
\bbe\label{eq:det}
\det\tau'(-1)=\det\tau(-1)^{[H_y:K]}=\det\tau(-1)^y.
\ee
\noindent
By Frobenius reciprocity
$$
\langle 1,\tau'\rangle+\langle\eta',\tau'\rangle=\langle
1_{H_y}\oplus\eta',\tau'\rangle=\langle
1_{H_y}\oplus\eta',\res^{H_y}_K\tau\rangle=\langle\ind^{H_y}_K(1_{H_y}\oplus\eta
'), \tau\rangle.
$$
As was mentioned above,
$\ind^{H_y}_K1_{H_y}=\bigoplus^{y-1}_{i=0}\chi_i$, where
$\chi_0,\ldots,\chi_{y-1}$ are all the unramified characters of
$K^{\times}$, satisfying $\chi_i(\varpi)^y=1$. Analogously,
$\ind^{H_y}_K\eta'=\bigoplus^{2y-1}_{i=y}\chi_i$, where
$\chi_y,\ldots,\chi_{2y-1}$ are all the unramified characters of
$K^{\times}$ satisfying $\chi_i^y(\varpi)=-1$ ($y\leq i\leq
2y-1$). Thus
$$
\ind^{H_y}_K(1_{H_y}\oplus\eta')=\bigoplus^{2y-1}_{i=0}\chi_i,
$$
where $\chi_0,\ldots,\chi_{2y-1}$ are all the unramified
characters of $K^{\times}$ satisfying $\chi_i(\varpi)^{2y}=1$, and
$$
\langle
1,\tau'\rangle+\langle\eta',\tau'\rangle=\sum_{i=0}^{2y-1}\langle\chi_i,\tau\rangle.
$$
Since $\tau$ has a real-valued character, for each $\chi_i$ of
order greater than 2, $\langle\chi_i,\tau\rangle$ will appear in
this sum twice, i.e.
\bbe\label{eq:4}
\langle 1,\tau'\rangle+\langle\eta',\tau'\rangle\equiv \langle
1,\tau\rangle+\langle\eta,\tau\rangle \qquad (\text{mod}\,\, 2).
\ee
Finally, by Frobenius reciprocity,
\bbe\label{eq:5}
\pair{\hat{\ga'},\tau'}=\pair{\ind^{H_{x}}_{H_y}(\phi\otimes\theta),\res^{H_y}_
K\tau}=
\pair{\ind^{H_{x}}_{K}(\phi\otimes\theta),\tau}=\pair{\hat\ga,\tau}.
\ee
Now formulas \eqref{eq:10} -- \eqref{eq:5} imply
\bbe\label{eq:lam}
W(\ga\otimes\tau)=\det\tau(-1)^y\cdot\varphi^{\dim\tau}\cdot
(-1)^{\pair{1,\tau}+\pair{\eta,\tau}+\pair{\hat\ga,\tau}}.
\ee
Applying \eqref{eq:lam} to $\mu_1,\ldots,\mu_a$ and substituting
the result into \eqref{eq:w} we get the statement of the
proposition.
\end{proof}

\section{}\label{ap:c}
\begin{lem}\label{l:ediv}
Let $G$ be a semi-abelian scheme over a field $K$. Then $G(\Kb)$
is a divisible abelian group.
\end{lem}
\begin{proof}
From \eqref{eq:uniform} we get the following exact sequence of
abelian groups:
$$
0\longrightarrow T(\Kb)\longrightarrow G(\Kb)\longrightarrow
A(\Kb) \longrightarrow 0.
$$
Clearly, $T(\Kb)$ is divisible, because $T(\Kb)\cong
(\Kb^{\times})^r$ as groups and $A(\Kb)$ is divisible, because $A$
is an abelian variety. Thus, $G(\Kb)$ is divisible as an abelian
group which is an extension of a divisible group by a divisible
group.
\end{proof}

In the next lemma we keep the notation of Subsection \ref{sub:loc gen}.
\begin{lem}\label{l:dimens}
Let $\sg'=(\sg,N)$ be the representation of $\wdk$ associated to the natural $l$-adic representation of $\gk$ on $V_l(M)^*$ and let $\sg'\cong \ga\oplus(\de\otimes\spin(2))$ for some representations
$\ga$ and $\de$ of $\wk$. Then $\dim\de=r$.
\end{lem}
\begin{proof}
Since $\dim\de=\rank N$ and for any finite extension $L\subset\Kb$ of $K$ we have $\res_{\cW'(\Kb/L)}\sg'=(\res^L_K\sg,N)$ (\cite{r2}, p. 130), by Lemma \ref{l:find} we can assume that $T$ splits over $K$ and $A$ has good reduction over $K$. 

We have the following exact sequence
of $\gk$-modules (\cite{ray}, p. 312):
\bbe\label{eq:al}
0\longrightarrow G(\Kb)_{l^n}\longrightarrow
M(\Kb)_{l^n}\rarab{\phi_{l^n}}Y(\Kb)/l^nY(\Kb)\rar 0.
\ee
Since $G(\Kb)$ is divisible by Lemma \ref{l:ediv}, sequence \eqref{eq:al} induces  an exact $\gk$-equivariant sequence of $l$-adic Tate modules:
$$
0\rar T_l(G)\rar T_l(M)\rar Z\rar 0,
$$
where $Z=\limproj (Y(\Kb)/l^nY(\Kb))$ with the maps being the natural
quotient maps. By tensoring the above sequence with $\bQ_l$ over $\bZ_l$ we
get the following exact $\gk$-equivariant sequence:
\bbe\label{eq:jkj}
0\rar V_l(G)\rar V_l(M)\rarab{\phi} Z\otimes_{\bZ_l}\bQ_l\rar 0.
\ee

Let $K^{tam}\subseteq\Kb$ be
the maximal tamely ramified extension of $K^{unr}$ and
$R=\Gal(\Kb/K^{tam})$. Then
$$
I/R\cong\prod_{\substack{s\text{ is a prime} \\ s\ne p}}\bZ_s.
$$
Let $t_l:I\rar\bZ_l$ denote the composition of the quotient map
onto $I/R$ with the projection onto $\bZ_l$.

Let $x_{l^n}\in M(\Kb)_{l^n}$, $i\in I$, and
$\phi_{l^n}(x_{l^n})=[y]$ for some $y\in Y(\Kb)$ and $\phi_{l^n}$ given by \eqref{eq:al}. Then a formula
on p. 314 in \cite{ray} yields:
\bbe\label{eq:x}
i(x_{l^n})=x_{l^n}+\nu_{l^n}(y\otimes t_{l^n}(i)),
\ee
where $t_{l^n}:I\rar\bZ/l^n\bZ$ is the composition of
$t_l:I\rar\bZ_l$ with the canonical projection onto the factor
$\bZ/l^n\bZ$;\quad $\nu_{l^n}:Y(\Kb)\otimes_{\bZ}(\bZ/l^n\bZ)\rar T(\Kb)_{l^n}$ \,is the following composition of
$\gk$-module homomorphisms:
$$
Y(\Kb)\otimes_{\bZ}(\bZ/l^n\bZ)\rar\Hom_{\bZ}(X(\Kb),\bZ)\otimes_{\bZ}(\bZ/l^n\bZ)\rarab{\sim}T(\Kb)_{l^n},
$$
where $X$ is the character group of $T$ and the first map is induced by the geometric monodromy
$$
\mu_0:Y\times X\rar\bZ;
$$
finally, $\nu_{l^n}(y\otimes t_{l^n}(i))\in T(\Kb)_{l^n}$ is
considered as an element of $M(\Kb)_{l^n}$ via the inclusions
$$
T(\Kb)_{l^n}\hookrightarrow G(\Kb)_{l^n}\hookrightarrow
M(\Kb)_{l^n}.
$$

For each $i\in I$ we have the following maps
$\al_n(i):Y(\Kb)/l^nY(\Kb)\rar T(\Kb)_{l^n}$ given by the
following composition:
$$
Y(\Kb)/l^nY(\Kb)\rarab{\psi_n(i)}Y(\Kb)\otimes_{\bZ}(\bZ/l^n\bZ)\rarab{\nu_{l^n}}
T(\Kb)_{l^n},
$$
where $\psi_n(i)([y])=y\otimes t_{l^n}(i)$, $y\in Y(\Kb)$. It is
easy to show that $\{\al_n(i)\}$ induce the homomorphism
$$
\al(i)=(\al_n(i)):Z\rar T_l(T),
$$
where $Z=\limproj(Y(\Kb)/l^nY(\Kb))$. By extending scalars to $\bQ_l$ we get
\bbe\label{eq:be}
\al'(i):Z\otimes_{\bZ_l}\bQ_l\rar V_l(T).
\ee

Let $\be_l:\gk\rar\GL(V_l(M))$ be the natural $l$-adic representation of
$\gk$ on $V_l(M)$. Then \eqref{eq:x} and \eqref{eq:be} imply:
$$
\be_l(i)=\id+\al'(i)\circ\phi,\quad i\in I,
$$
where $\id:V_l(M)\rar V_l(M)$ is the identity map and $\phi$ is given by \eqref{eq:jkj}. On the other
hand,
$$
\be_l(i)=\exp(t_l(i)R_l),
$$
where $i$ is in some open subgroup $J$ of $I$ and $R_l$ is a
nilpotent endomorphism on $V_l(M)$ (\cite{r2}, Prop. on p. 131). Since $\sg'=(\sg,N)$ is the representation of $\wdk$ associated to the representation $\sg_l':\gk\rar\GL(V_l(M)^*)$, it follows that $N$ is obtained from $-R_l^t$ by extending scalars via a field embedding $\imath:\bQ_l\hookrightarrow \bC$.
Thus, $\rank R_l=\rank N=\dim\de$ and $R_l^2=N^2=0$ by assumption.
Thus,
$$
\al'(i)\circ\phi=t_l(i)R_l,\quad i\in J,
$$
and it is enough to show that there exists $i_0\in I$ such that
$\al'(i_0)$ is surjective.

For each $n$ let $\xi_n$ be a generator of $\bZ/l^n\bZ$ such that
$l\cdot\xi_{n+1}=\xi_n$ so $\xi=(\xi_n)\in\bZ_l$. Since
$t_l:I\rar\bZ_l$ is surjective, there exists $i_0\in I$ such that
$t_l(i_0)=\xi$. It implies that $\psi_n(i_0)$ is an isomorphism
for each $n$, hence it is enough to show that the map
$$
\nu':\limproj(Y(\Kb)\otimes_{\bZ}(\bZ/l^n\bZ))\otimes_{\bZ_l}\bQ_l\rar
V_l(T)
$$
induced by $(\nu_{l^n})$ is surjective. Since $\mu_0$ is
nondegenerate (\cite{ch}, p. 52, Remark 6.3), we have the
following exact sequence:
$$
0\rar Y(\Kb)\rarab{g}\Hom_{\bZ}(X(\Kb),\bZ)\rar B\rar 0,
$$
where $g(y)=\mu_0(y,\cdot)$ and $B$ is finite, since $Y(\Kb)$
and $\Hom_{\bZ}(X(\Kb),\bZ)$ are free abelian groups of the same
rank $r$. Applying the functor $(-)\otimes_{\bZ}(\bZ/l^n\bZ)$ to
the above sequence, we get:
$$
Y(\Kb)\otimes_{\bZ}(\bZ/l^n\bZ)\rarab{\nu_{l^n}} T(\Kb)_{l^n}\rar
B\otimes_{\bZ}(\bZ/l^n\bZ)\rar 0,
$$
hence
$$
0\rar\im\nu_{l^n}\rar T(\Kb)_{l^n}\rar B/l^nB\rar 0.
$$
Since $Y(\Kb)\otimes_{\bZ}(\bZ/l^n\bZ)$ is a finite group,
$\im\nu_{l^n}$ is a finite group, hence $\{\im\nu_{l^n}\}$
satisfies the Mittag-Leffler condition and we have the following
exact sequence:
\bbe\label{eq:ga}
0\rar\limproj(\im\nu_{l^n})\rar T_l(T)\rar\limproj(B/l^nB)\rar 0.
\ee
Here $\limproj(\im\nu_{l^n})\cong\im\nu$, where
$$\nu=(\nu_{l^n}):\limproj(Y(\Kb)\otimes_{\bZ}(\bZ/l^n\bZ))\rar T_l(T).$$
Indeed, let $S_n=Y(\Kb)\otimes_{\bZ}(\bZ/l^n\bZ)$, then we have an
exact sequence
$$
0\rar\ker\nu_{l^n}\rar S_n\rar\im\nu_{l^n}\rar 0,
$$
where the maps from $S_n$ to $\im\nu_{l^n}$ are induced by $\nu_{l^n}$.
Since $\ker\nu_{l^n}$ is finite for each $n$, $\{\ker\nu_{l^n}\}$
satisfies the Mittag-Leffler condition, hence there is the
following exact sequence
\bbe\label{eq:ex}
0\rar\limproj(\ker\nu_{l^n})\rar\limproj
S_n\rar\limproj(\im\nu_{l^n})\rar 0.
\ee
On the other hand, from \eqref{eq:ga} we have
$$
\limproj(\im\nu_{l^n})\hookrightarrow T_l(T),
$$
which together with \eqref{eq:ex} implies
$\limproj(\im\nu_{l^n})\cong\im\nu$.

Thus, applying the exact functor $(-)\otimes_{\bZ_l}\bQ_l$ to
\eqref{eq:ga} and taking into account that $B$ is finite, we get:
$$
0\rar(\im\nu)\otimes_{\bZ_l}\bQ_l\rar V_l(T)\rar 0,
$$
which implies that
$$
\im\nu'\cong (\im\nu)\otimes_{\bZ_l}\bQ_l\cong V_l(T),
$$
hence $\nu'$ is surjective.
\end{proof}

\section{}\label{ap:a}
In this appendix we keep the notation of Subsection \ref{sub:loc gen}. Let
$\rho'=(\rho,S)$ be the representation of $\wdk$ associated to the
natural $l$-adic representation  of $\gk$ on $V_l(T(\Kb)/\La)^*$,
where $\La=T(\Kb)\bigcap Y(\Kb)$ is a free discrete subgroup of
$T(\Kb)$ of rank $s$ ($s\leq r$). It is known that there is a
finite Galois extension $L\subset\Kb$ of $K$ such that
$\Gal(\Kb/L)$ acts trivially on $Y(\Kb)$, hence $Y(\Kb)$ can be
considered as a $\Gal(L/K)$-module. Let
$$\chi:\Gal(L/K)\rar\GL_r(\bZ)$$
denote the corresponding representation. Thus, from
\eqref{eq:uniform} we have the following exact sequence of
$\Gal(L/K)$-modules:
$$
0\longrightarrow \La\otimes_{\bZ}\bC\longrightarrow
Y(\Kb)\otimes_{\bZ}\bC\longrightarrow
B\otimes_{\bZ}\bC\longrightarrow 0,
$$
where $B=f(Y(\Kb))$. Let $\chi_1:\Gal(L/K)\longrightarrow
\GL_s(\bZ)$ denote the representation of $\Gal(L/K)$ on $\La$ and
$\chi_2:\Gal(L/K)\longrightarrow \GL_{r-s}(\bZ)$ denote the
representation of $\Gal(L/K)$ on $B\otimes_{\bZ}\bC$. Then
$$
\chi\cong\chi_1\oplus\chi_2.
$$

\begin{prop}\label{p:nonsplit tor}
$$
\rho'\cong(\chi_2\otimes\om^{-1})\oplus(\chi_1\otimes\om^{-1}\otimes
\spin(2)).
$$
\end{prop}

\begin{cor}\label{cor:root}
Let $\tau$ be a representation of $\gk$ with real-valued
character. If $r=s$, then $\chi=\chi _1$ and we have:
$$
W(\rho'\tens\tau)=\det\tau(-1)^r\cdot\det\chi(-1)^{\dim\tau}\cdot
(-1)^{\langle\chi,\tau\rangle}.
$$
\end{cor}
\begin{proof}
The statement is a consequence of Proposition \ref{p:nonsplit tor}
and Proposition 6 on p. 327 in \cite{r1} (cf. \cite{r1}, p. 329,
Thm. 2(ii)).
\end{proof}

\begin{proof}[Proof of Proposition \ref{p:nonsplit tor}]
Let $\Ga=T(\Kb)/\La$. We have the following exact sequence of
$\gk$-modules:
\bbe
0 \rar\La \rar T(\Kb) \rar \Ga \rar 0.
\ee
Since $\La\cong\bZ^s$ and $T(\Kb)$ is a divisible group, this
sequence induces the following exact $\gk$-equivariant sequence of
$l$-adic Tate modules:
\bbe\label{eq:zvezda}
0 \rar T_l(T) \rar T_l(\Ga) \rar \chi_1\tens\bZ_l^s \rar 0,
\ee
where $T_l(T)$ denotes $T_l(T(\Kb))$. 
Let $L\subset\Kb$ be a finite Galois extension of $K$ over which
$T$ splits. Since $T_l(T)$ is a free $\bZ_l$-module of rank $r$ it follows
from \eqref{eq:zvezda} that $T_l(\Ga)$ is a free $\bZ_l$-module of
rank $s+r$, hence by Proposition \ref{p:splittor} below we have 
\bbe\label{eq:h}
\res^L_K\rho'\cong(\om_L^{-1})^{\oplus(r-s)}\oplus(\om_L^{-1}\tens\spin(2))^{\oplus
s},
\ee
where $\om _L=\res^L_K\om$. The rest of the proof is similar to the proof of Proposition \ref{pr:exact_sequence}. Since $\res_{\cW'(\Kb/L)}\rho'=(\res^L_K\rho,S)$ and $\res^L_K\rho$ is semisimple by \eqref{eq:h}, $\rho'$ is admissible by Lemma \ref{l:find}. Hence it
has the following form:
\bbe\label{eq:ra}
\rho'\cong\bigoplus_{i=1}^t\al_t\otimes\spin(n_t),
\ee
where each $\al_i$ is a representation of $\wk$ and each $n_i$ is a positive integer (\cite{r2}, p. 133, Cor. 2). Also, it follows from \eqref{eq:h} that $S^2=0$ and $\rank S=s$. Thus, each $n_i$ in \eqref{eq:ra} is 1 or 2 and
\bbe\label{eq:1}
\rho'\cong\al\oplus(\be\otimes\spin(2)),
\ee
where $\al$ is a representation of $\wk$ of dimension $r-s$ and
$\be$ is a  representation of $\wk$ of dimension $s$.

Applying the exact functor
$(-)\tens_{\bZ_l}\bQ_l$ to \eqref{eq:zvezda} and taking duals
afterwards, we get
\bbe\label{eq:B}
0 \rar \chi_1\tens\bQ_l^s\rar V_l(\Ga)^* \rar V_l(T)^* \rar 0,
\ee
where $V_l(T)=V_l(T(\Kb))$ and $\chi_1\cong(\chi_1)^*$, since
$\chi_1$ is a representation of finite image, realizable over
$\bZ$. 
Sequence \eqref{eq:B} induces an exact sequence of corresponding representations of $\wdk$, i.e.
\bbe\label{eq:B1}
0\longrightarrow (\chi_1\otimes\bQ_l^s)\otimes_{\imath}\bC\longrightarrow
V_l(\Ga)^*\otimes_{\imath}\bC\rar V_l(T)^*\otimes_{\imath}\bC
\longrightarrow 0
\ee
is an exact sequence of $\wdk$-modules, where $\chi_1$ is the representation of $\wdk$ on $(\chi_1\otimes\bQ_l^s)\otimes_{\imath}\bC$, $\rho'=(\rho,S)$ is the representation of $\wdk$ on $V_l(\Ga)^*\otimes_{\imath}\bC$, and by Lemma \ref{l:ator}, $\chi\otimes\om^{-1}$ is the representation of $\wdk$ on $V_l(T)^*\otimes_{\imath}\bC$. Since $\rho$ is semisimple, the exact sequence \eqref{eq:B1} of $\wk$-modules splits, i.e.
$$
\rho\cong\chi_1\oplus(\chi\otimes\om^{-1}).
$$
On the other hand, from \eqref{eq:1} we have:
$$
\rho\cong\al\oplus\be\oplus(\be\otimes\om).
$$
Thus, combining the last two congruences, we get
\bbe\label{eq:tri}
\al\oplus\be\oplus(\be\otimes\om)\cong\chi_1\oplus(\chi\otimes\om^{-1}).
\ee

We claim that $\be\otimes\om$ is
isomorphic to a subrepresentation of $\chi_1$. Suppose there is an irreducible component $\be_0$ of $\be$ such that $\be_0\otimes\om$ is
isomorphic to a subrepresentation of $\chi\otimes\om^{-1}$, i.e.
\bbe\label{eq:odinn}
\be_0\otimes\om\cong x\otimes\om^{-1}
\ee
for some irreducible component $x$ of $\chi$.
It follows from \eqref{eq:tri} that $\be_0$ is isomorphic to a subrepresentation of $\chi\otimes\om^{-1}$ or $\chi_1$, which is impossible, because
$x$, $\chi$, and $\chi_1$ have finite images, whereas $\om$ does
not. Indeed, suppose $\be_0\cong y\otimes\om^{-1}$ or $\be_0\cong
z$, where $y$ is an irreducible component of $\chi$ and $z$ is an
irreducible component of $\chi_1$. From \eqref{eq:odinn} we get
$$
\be_0\cong x\otimes\om^{-2},
$$
hence
$$
x\otimes\om^{-2}\cong y\otimes\om^{-1}\quad\text{or}\quad
x\otimes\om^{-2}\cong z.
$$
By taking determinants of both sides in each case, we get
\bbe\label{eq:opo}
\frac{\det x}{\det y}=\om^{-m+2k} \quad\mbox{or}\quad \frac{\det
x}{\det z}=\om^{2k},
\ee
where $m=\dim y$, $k=\dim x$. Since $x$, $y$, and $z$ have finite
images (as being subrepresentations of $\chi$ or $\chi_1$) and
$\om$ has infinite image, \eqref{eq:opo} gives a contradiction.

Thus, $\be\otimes\om$ is isomorphic to a subrepresentation of $\chi_1$.
Since $\be\otimes\om$ and $\chi_1$ have the same dimension $s$, we have
$\be\otimes\om\cong\chi_1$, hence $\be\cong\chi_1\otimes\om^{-1}$.
By the uniqueness of decomposition of a semisimple module into simple modules, we conclude from \eqref{eq:tri} that $\al\cong\chi_2\otimes\om^{-1}$.
\end{proof}

\begin{prop}\label{p:splittor}
Let $\La\subset\ktt$ be a free
discrete subgroup of rank $s$ $(s\leq r)$ and denote $\kbtt/\La$
by $\Ga$. Let $\rho'=(\rho,S)$ be the representation of
$\cW'(\Kb/K)$ associated to the $l$-adic representation of $\gk$
on $V_l(\Ga)^*$. Let $T_l(\Ga)$ be a free $\bZ_l$-module of rank
$s+r$. Then
$$
\rho'\cong(\om^{-1})^{\oplus(r-s)}\oplus(\om^{-1}\otimes
\spin(2))^{\oplus s}.
$$
\end{prop}

\begin{proof} Let $p_1,\dotsc,p_s\in\La$ be a basis of $\La$, satisfying the
assertion of Lemma \ref{l:basis} below. First, let us choose a
$\bQ_l$-basis for $V_l=V_l(\Ga)$. Let
$f_1=(f_1(n)),\dotsc,f_s=(f_s(n))\in T_l(\Ga)$, where
$f_1(n),\dotsc,f_s(n)$ as elements of $\kbtt$ have the following
form:
\begin{eqnarray*}
f_1(n)^{l^n}  =  p_1, & \dotsc &, f_s(n)^{l^n} = p_s  \quad\text{and} \\
f_1(n+1)^l  =  f_1(n), & \dotsc &, f_s(n+1)^l = f_s(n).
\end{eqnarray*}
Let $\xi=(\xi(n))$, where each $\xi(n)\in\kbt$ is a primitive
$l^n$-th root of unity and $\xi(n+1)^l=\xi(n)$. Let
$f_{s+1}=(f_{s+1}(n)),\dotsc,f_{s+r}=(f_{s+r}(n))\in T_l(\Ga)$,
where $f_{s+1}(n),\dotsc,f_{s+r}(n)$ as elements of $\kbtt$
satisfy the following properties:
\begin{eqnarray*}
f_{s+1}(n)&=&(\xi(n),1,\dotsc,1),\\
f_{s+2}(n)&=&(1,\xi(n),\dotsc,1), \\
\dotsc \\
f_{s+r}(n)&=&(1,1,\dotsc,\xi(n)).
\end{eqnarray*}
Then $f_1,\dotsc,f_{s+r}$ is a basis of $V_l$. Indeed, it is easy
to check that $f_1,\dotsc,f_{s+r}$ are linearly independent over
$\bZ_l$. Since $T_l(\Ga)$ is a free $\bZ_l$-module of rank $s+r$,
it follows that $f_1,\dotsc,f_{s+r}$ is a basis of $V_l$.

Let $\rho_l:\gk\rar\GL(V_l)$ be the $l$-adic representation
associated to the $\gk$-module $V_l$. Then the matrix
representation of $\rho_l$ with respect to the basis
$f_1,\dotsc,f_{s+r}$ has the following form:
\bbe\label{eq:tor}
{\rho_l(\Phi)}= \left(\begin{array}{cc}
E_s & 0 \\
* & q^{-1}\cdot E_r\\
\end{array} \right), \quad
{\rho_l(i)}= \left(\begin{array}{cc}
E_s & 0\\
B(i) & E_r\\
\end{array} \right),
\ee
\noindent
where $E_r$ and $E_s$ are the identity matrices, $i\in I$, and
$B(i)\in\text{Mat}_{r\times s}(\bQ_l)$. It is known that there
exists a nilpotent endomorphism $S_l$ of $V_l^*$ such that $S$ is
obtained from $S_l$ by extending of scalars via a field embedding
$\imath:\bQ_l\hookrightarrow\bC$; moreover, $S_l$ is a unique
nilpotent endomorphism such that
\bbe\label{eq:nil}
\rho_l^*(i)=\text{exp}(t_l(i)S_l),
\ee
\noindent
where $t_l:I\rar\bQ_l$ is a nontrivial continuous homomorphism and
$i$ belongs to an open subgroup of $I$. Furthermore, for any
$g\in\wk$ we have
\bbe\label{eq:weilgr}
\rho(g)=\rho_l^*(g)\text{exp}(-t_l(i)S_l),
\ee
\noindent
where $\rho_l^*(g)\text{exp}(-t_l(i)S_l)$ is considered as an
element of $\GL(V_l^*\otimes_{\imath}\bC)$ via $\imath$
(\cite{r2}, p. 131, Prop.(i), (ii)). Formula \eqref{eq:tor} for
$\rho_l(\Phi)$ implies that, considered as a matrix over $\bC$ via
$\imath$, it is diagonalizable. It follows from Formula
\eqref{eq:weilgr} that $\rho(\Phi)$ is diagonalizable, hence
$\rho$ is semisimple and $\rho'$ is admissible by Lemma \ref{l:find}.

Let $\vp$ be a uniformizer of $K$ and let $\Pi=(\vp(n))$, where
each $\vp(n)\in\kbt$ has the following property:
$$
\vp(n)^{l^n}=\vp \quad \text{and} \quad \vp(n+1)^l=\vp(n).
$$
\noindent
There exists $i_0\in I$ such that
$$
i_0(\Pi)=(i_0(\vp(n)))=(\vp(n)\xi(n)^{\al(n)})=\xi^{\al}\Pi,
$$
\noindent
where $\al=(\al(n))\in\bZ_l$. By Lemma \ref{l:simple} below
$\al\ne 0$.
\begin{lem}\label{l:simple}
Let $g\in\cO$ and let $g_n\in\Kb$ denote a root of $x^{l^n}-g=0$.
Then $i(g_n)=g_n$ for any $i\in I$ and $n\in\bN$ if and only if
$g\in\cO^{\times}$.
\end{lem}
\begin{proof}
Clearly, $i(g_n)=g_n$ for any $i\in I$ and $n\in\bN$ if and only
if $K(g_n)$ is unramified over $K$ for any $n$.

Let $g\in\cO^{\times}$. Then the assertion follows from the fact
that $x^{l^n}-g$, considered as a polynomial in $k[x]$, has no
multiple root (\cite{l}, p. 48, Prop. 7).

Conversely, since $g_n^{l^n}=g$, the valuation of $g$ in $K(g_n)$
is divisible by $l^n$. Since $K(g_n)$ is unramified over $K$ for
any $n$, the valuation $v_g$ of $g$ in $K$ coincides with the
valuation of $g$ in $K(g_n)$, hence $v_g$ is divisible by $l^n$
for any n, which implies that $v_g$ must be zero and
$g\in\cO^{\times}$.
\end{proof}

Let $p_k=(p_{kj})$, $p_{kj}\in K^{\times}$, $1\leq k\leq s$,
$1\leq j\leq r$. It follows from Lemma \ref{l:basis} below that without
loss of generality we can assume that $p_{kk}\not\in\ot$ for any
$k$ and that $p_{kj}\in\ot$ whenever $k>j$, $1\leq j\leq r$. Thus
there exist $u_k\in\ot$ and $m_k\in\bZ^{\times}$ such that
$p_{kk}=u_k\cdot\vp^{m_k}$. Let $(u_k(n))$ be a sequence in $\kbt$
such that
$$
u_k(n)^{l^n}=u_k \quad \text{and} \quad u_k(n+1)^l=u_k(n).
$$
For $f_k(n)\in\kbtt$ write $f_k(n)=(f_{kj}(n))$, where
$f_{kj}(n)\in\kbt$, $1\leq k\leq s$, $1\leq j\leq r$.
\noindent
Then as $f_{kk}(n)$ we can take $u_k(n)\cdot\vp(n)^{m_k}$. For
$i_0$ we have
$$
i_0(f_1)=(i_0(f_1(n)))=(i_0(f_{11}(n)),i_0(f_{12}(n)),\dotsc,i_0(f_{1r}(n))),
$$
\noindent
where by Lemma \ref{l:simple} we have:
$$
i_0(f_{11}(n))=u_1(n)\cdot\vp(n)^{m_1}\cdot\xi(n)^{\al(n)
m_1}=f_{11}(n)\cdot\xi(n)^{\al(n) m_1}.
$$
Analogously, using Lemma \ref{l:simple}, we get the following
formulas:
\begin{eqnarray*}
i_0(f_1) & = & f_1\cdot f_{s+1}^{\al m_1}\cdot f_{s+2}^{a_2}\dotsb f_{s+r}^{a_r},\\
i_0(f_2) & = & f_2\cdot f_{s+2}^{\al m_2}\cdot f_{s+3}^{b_3}\dotsb f_{s+r}^{b_r},\\
\dotsc \\
i_0(f_s) & = & f_s\cdot f_{2s}^{\al m_s}\cdot
f_{2s+1}^{c_{s+1}}\dotsb f_{s+r}^{c_r}
\end{eqnarray*}
for some $a_i,b_j,\ldots,c_k\in\bZ_l$. This implies that $B(i_0)$
has the following form:
\[
{B(i_0)}= \left(\begin{array}{cccc}
\al m_1 & 0 & \dotsc & 0\\
* & \al m_2 & \dotsc & 0\\
\vdots & \vdots & \ddots & \vdots\\
* & * & \dotsc & \al m_s \\
\vdots & \vdots & \ddots & \vdots
\end{array} \right),
\]
\noindent
where $\al,m_1,\dotsc,m_s\in\bZ_l^{\times}$, hence $\rank
B(i_0)=s$.

Since $\rho'$ is admissible,
\bbe\label{eq:decomp}
\rho'\cong\bigoplus_{j=1}^{k}\pi_{j}\otimes\spin(n_j),
\ee
\noindent
where $\pi_1,\dotsc,\pi_k$ are representations of $\wk$
(\cite{r2}, p. 133, Cor. 2 ). Since
$\rho^*_l(i)=\text{exp}(t_l(i)S_l)$ by \eqref{eq:nil} and
$(\rho_l(i)-E_{s+r})^2=0$ from \eqref{eq:tor}, it follows that
$(S_l)^2=0$, i.e. each $n_j$ in \eqref{eq:decomp} is $1$ or $2$,
hence
$$
\rho'\cong\al\oplus(\be\otimes\spin(2)),
$$
where $\al$ and $\be$ are representations of $\wk$. Since $\rank
B(i_0)=s$, the equation \eqref{eq:nil} implies that $\rank S=\rank
S_l=s$, hence $\dim\be=s$ and $\dim\al=r-s$.

Let us prove now that $\al\cong\bigoplus_{r-s}\om^{-1}$ and
$\be\cong\bigoplus_s\om^{-1}$. It can be easily verified that
\eqref{eq:tor} -- \eqref{eq:weilgr} imply
\bbe\label{eq:star}
\rho(g)=\left(\begin{array}{cc}
E_s & * \\
0 & (\om^{-1})^{\oplus r}\end{array}\right),\quad g\in\wk,
\ee
hence there is a complete flag of subrepresentations
$$(0)\neq
W_1\subset\dotsb\subset W_{s+r}=V_l^*\otimes_{\imath}\bC$$ of
$\rho$. Since $\rho$ is semisimple, it implies that $\rho$ is a
direct sum of one-dimensional subrepresentations, hence from
\eqref{eq:star}
$$
\al\cong\bigoplus_{r-s}\om^{-1} \quad\text{and}\quad
\be\cong\bigoplus_s\om^{-1}.
$$
\end{proof}

\begin{lem}\label{l:basis}
Let $K$ be a non-Archimedean local field of characteristic zero
with ring of integers $\cO$. Let $\La\subset(K^{\times})^r$ be a
free discrete subgroup of rank $s$ ($s\leq r$). There exist a
basis $p_1,\ldots,p_s$ of $\La$ and natural numbers
$n_1,\ldots,n_s$ ($1\leq n_1<n_2<\dotsb <n_s\leq r$) with the
following property: if $p_k=(p_{kj})$, $1\leq k\leq s$, $1\leq j
\leq r$, and $p_{kj}\in\kt$, then $p_{in_i}\not\in\ot$ for any $i$
and $p_{ln_i}\in\ot$ whenever $l>i$.
\end{lem}

\begin{proof}
First, note that $\ol=\{1\}$. Indeed, if $x\in\ol$ and $x\neq 1$,
then $(x^n)$ is an infinite sequence in $\ol$, hence it has a
limit point, because $\ott$ is compact, which contradicts the
assumption that $\La$ is discrete.

Let $\varpi$ be a uniformizer of $K$. The map
$\ot\times\bZ\longrightarrow\kt$ given by
$$
(u,n)\mapsto u\varpi^n
$$
is an isomorphism of topological groups. For every positive
integer $r$ it induces an isomorphism $\ktt\cong\ott\times\bZ^r$.
Let $\pi:\ktt\rar\bZ^r$ be the projection  onto $\bZ^r$ and
$t_1,\ldots,t_s$ be a basis of $\La$. Since $\ol=\{1\}$,
$\pi(t_1),\ldots,\pi(t_s)$ form a basis of $\pi(\La)$. Indeed,
otherwise, there exist $m_1,\dotsc,m_s\in\bZ$, not all of which
are zeros, such that
$$
m_1\pi(t_1)+\dotsb+m_s\pi(t_s)=0.
$$
\noindent
Then $1\neq t_1^{m_1}\dotsm t_s^{m_s}\in\ol$. Thus,
$\pi(\La)\subseteq\bZ^r$ is a subgroup of rank $s$ and it is
enough to prove the following sublemma:

\begin{sublem}\label{sublemmma}
Let $G\subseteq\bZ^r$ be a subgroup of rank $s$ $(s\leq r)$. There
exist a basis $g_1,\ldots,g_s$ of $G$ and natural numbers
$n_1,\ldots,n_s$ $(1\leq n_1<n_2<\dotsb <n_s\leq r)$ with the
following property: if $g_k=(g_{kj})$, $1\leq k\leq s$, $1\leq j
\leq r$, and $g_{kj}\in\bZ$, then $g_{in_i}\ne 0$ for any $i$ and
$g_{ln_i}=0$ whenever $l>i$.
\end{sublem}

Indeed, if we assume Sublemma \ref{sublemmma}, then there is a
basis $g_1,\ldots,g_s$ of $\pi(\La)$ with the property described in Sublemma \ref{sublemmma}. Since $\pi(t_1),\ldots,\pi(t_s)$ is a basis of $\pi(\La)$, there is a matrix $D=(d_{ij})\in\GL_s(\bZ)$ such that 
$$g_i=\sum_j
d_{ij}\pi(t_j),\quad 1\leq i\leq s.
$$
Then $p_i=\prod_j t_j^{d_{ij}}$, $1\leq
i\leq s$, will be a basis of $\La$ with the required property.
\end{proof}

\begin{proof}[Proof of Sublemma \ref{sublemmma}]
Suppose $r=s$. We will prove the sublemma in this case by
induction on $r$. Clearly, it holds when $r=1$. Let $r$ be
arbitrary and $e_1,\dotsc,e_r$ be the standard basis of $\bZ^r$.
There is $k\in \bZ^{\times}$ such that $G\cap(e_r)=(ke_r)$,
because $me_r\in G$, where $m=|B|$ and $B=\bZ^r/G$. Then
$G/(ke_r)\subseteq \bZ^{r-1}$ is a subgroup of rank $r-1$. By
induction, there exist $g_1,\ldots,g_{r-1}\in G$ such that in
$G/(ke_r)$ we have:
$$
\bar{g_i}=\sum_{j=1}^{r-1} a_{ij}e_j, \quad 1\leq i\leq r-1,
$$
for some $a_{ij}\in\bZ$ such that $a_{ii}\ne 0$ for any $i$ and
$a_{ij}=0$ whenever $i>j$. Then $g_1,\ldots,g_{r-1},g_r=ke_r$ will
be a basis of $G$ with the required property.

Suppose now that $s\ne r$. Let $q_1,\ldots,q_s$ be a basis of $G$.
Then $q_i=\sum_j b_{ij}e_j$, where
$B=(b_{ij})\in\text{Mat}_{s\times r}(\bZ)$. Since $q_1,\ldots,q_s$
is a basis, $\rank B=s$, i.e. there exists an $s\times s$-submatrix $B_0$ of
$B$ such that $\det B_0\ne 0$. Let $B_0$ have the following form:
\begin{displaymath}
{B_0}=\left( \begin{array}{cccc}
b_{1n_1} & b_{1n_2} & \ldots & b_{1n_s} \\
b_{2n_1} & b_{2n_2} & \ldots & b_{2n_s} \\
\vdots & \vdots & \ddots & \vdots \\
b_{sn_1} & b_{sn_2} & \ldots & b_{sn_s}\end{array}\right) .
\end{displaymath}
Then $p_i=\sum_j b_{in_j}e_{n_j}$, $1\leq i\leq s$, are linearly
independent, hence generate a free subgroup $H$ of rank $s$ in
$\bZ e_{n_1}\oplus\dotsb\oplus\bZ e_{n_s}$. By the case $r=s$ above there is a
matrix  $C\in\text{GL}_s(\bZ)$ such that
$$
\sum_i c_{ki}p_i=\sum_j h_{kn_j}e_{n_j},
$$
where $h_{kn_j}\in\bZ$, $h_{in_i}\ne 0$ for any $i$ and
$h_{ln_i}=0$ whenever $l>i$. Then $g_k=\sum_i c_{ki}q_i$, $1\leq
k\leq s$, will be a basis of $G$ with the required property.
\end{proof}

\section{}\label{ap:d}

We keep the notation of Subsection \ref{sub:loc gen} except that $K$ does not have to be of characteristic zero and $\Kb$ denotes a separable algebraic closure of $K$.
As in Appendix \ref{sec:b} (see Lemma \ref{l:sympl}) if
$\la:D\longrightarrow \GL(U)$ is a representation of a group $D$,
then by $\check\la:D\longrightarrow \GL(\check U)$ we denote the
representation of $D$ on $\check U$, where $\check U$ is a
$\bC[D]$-module with the underlying $D$-module $U^*$ and
multiplication by constants defined as follows:
$$
a\cdot\phi=\overline{a}\phi,\quad a\in\bC,\,\phi\in U^*.
$$
We say that $U$ is {\it unitary} if $U$ admits a nondegenerate
invariant hermitian form (not necessarily positive definite).

\begin{prop}\label{pr:lok}
Let $\sg'$ be an admissible representation of $\wdk$ written in
the following form:
$$
\sg'\cong\bigoplus_{i=1}^{k}\pi_i\otimes\spin(n_i),
$$
where each $\pi_i$ is a representation of $\wk$ and $n_i\ne n_j$ whenever $i\ne j$ $($\cite{r2}, p. $133$, Cor. $2)$. If $\sg'$ is unitary,
orthogonal, or symplectic with respect to a corresponding
invariant nondegenerate form $\langle\cdot,\cdot\rangle$ then
each $\pi_i\otimes\spin(n_i)$ is unitary, orthogonal, or
symplectic respectively with respect to the restriction of
$\langle\cdot,\cdot\rangle$.
\end{prop}
\begin{proof}
Let $U$ be a representation space of $\sg'$ and $U_i$ a
representation space of $\pi_i\otimes\spin(n_i)$, $1\leq i\leq k$, so that $U=\bigoplus_{i=1}^k U_i$.
Let $\langle\cdot\, ,\cdot\rangle$ be a nondegenerate invariant
form on $U$ and let $\tilde U$ be a $\wdk$-module over $\bC$ such that
$\tilde U=U^*$ if $\langle\cdot\, ,\cdot\rangle$ is bilinear and
$\tilde U=\check U$ if $\langle\cdot\, ,\cdot\rangle$ is
sesquilinear. Let $\phi:U\longrightarrow\tilde  U$ be the
isomorphism of $\wdk$-modules induced by $\langle\cdot\,
,\cdot\rangle$, and let $\psi:\tilde
U={(\bigoplus_{i=1}^k U_i)}^{\sim}\longrightarrow\tilde{U_1}\oplus\dotsb\oplus\tilde{U_k}$ denote
the usual isomorphism. It is easy to show that for any $n$ we have
${\spin(n)}^{\sim}\cong\om^{-(n-1)}\otimes\spin(n)$, hence
$\tilde{U_i}\cong V_i$, where $V_i$ denotes a representation space
of $\tilde{\pi_i}\otimes\om^{-(n_i-1)}\otimes\spin(n_i)$. Denote
by $\la:\tilde{U_1}\oplus\dotsb\oplus\tilde{U_k}\longrightarrow
V_1\oplus\dotsb\oplus V_k$ the corresponding isomorphism. For each
$i$ let $\rho_i:U_i\longrightarrow V_i$ be defined by the
following diagram:
$$
\xymatrix{
U \ar[r]^-{\la\circ\psi\circ\phi}   &   V_1\oplus\dotsb\oplus V_k
\ar[d]^{\pi_i} \\
U_i \ar@{^{(}->}[u] \ar[r]^{\rho_i} &  V_i }
$$
where $\pi _i$ is the projection onto $i$-th factor. To prove that
$\langle\cdot\, ,\cdot\rangle\vert_{U_i}$ is nondegenerate for
each $i$ is equivalent to proving that $\rho_i$ is an isomorphism
for each $i$, which follows from Lemma \ref{l:appc} below.
\end{proof}

\begin{lem}\label{l:appc}
Let $\al'=(\al,N)$ and $\be'=(\be,M)$ be two isomorphic admissible
representations of $\wdk$ written in the following forms:
$$
\al'\cong\bigoplus_{i=1}^{k}\al_i\otimes\spin(n_i)\quad\text{and}\quad\be'\cong\bigoplus_{i=1}^{k}\be_i\otimes\spin(n_i),
$$
where each $\al_i$ and $\be_i$ is a representation of $\wk$, $n_i\ne n_j$ whenever $i\ne j$. Let $U$ $($resp. $V)$ be a
representation space of $\al'$ $($resp. $\be')$ and for each $i$
let $U_i$ $($resp. $V_i)$ be a representation space of
$\al_i\otimes\spin(n_i)$ $($resp. $\be_i\otimes\spin(n_i))$. Let
$\phi:U\longrightarrow V$ be an isomorphism of $\wdk$-modules and
$\psi_i:U_i\longrightarrow V_i$, $1\leq i\leq k$, defined by the
following diagram:
$$
\xymatrix{
U \ar[r]^-{\phi}   &   V
\ar[d]^{\pi_i} \\
U_i \ar@{^{(}->}[u] \ar[r]^{\psi_i} &  V_i }
$$
where $\pi_i$ is the projection onto $i$-th factor. Then each
$\psi_i$ is an isomorphism of $\wdk$-modules.
\end{lem}

\begin{proof}
We will prove the lemma by induction on $k$. Clearly, it holds
when $k=1$. Let $k$ be arbitrary and let $e_0,\ldots,e_{n_k-1}$ be the
standard basis of $\bC^{n_k}$. Without loss of generality we can
assume that $n_k>n_i$ for any $i$. Let $U_k^{\circ}$ be a
representation space of $\al_k$, so $U_k=\bigoplus
_{j=0}^{n_k-1}(U_k^{\circ}\otimes e_j)$. Then
\begin{eqnarray}
U & = & \ker N^{n_k-1}\oplus (U_k^{\circ}\otimes e_0)\quad\text{and} \nonumber \\
N^jU & = & N^j(\ker N^{n_k-1})\oplus (U_k^{\circ}\otimes
e_j),\label{eq:tar}\quad 0\leq j\leq n_k-1.
\end{eqnarray}
Since $\phi$ is an isomorphism of $\wdk$-modules, we have from
\eqref{eq:tar}:
\bbe\label{eq:(1)}
M^jV=M^j(\ker M^{n_k-1})\oplus \phi(U_k^{\circ}\otimes e_j),\quad
0\leq j\leq n_k-1.
\ee
On the other hand, \eqref{eq:tar} holds in $V$, i.e.
\bbe\label{eq:(2)}
M^jV=M^j(\ker M^{n_k-1})\oplus (V_k^{\circ}\otimes e_j),\quad
0\leq j\leq n_k-1,
\ee
where $V_k^{\circ}$ denotes a representation space of $\be_k$. We
have the following filtration of $V$:
\begin{eqnarray}
V\supseteq\ker M^{n_k-1}\supseteq MV\supseteq & M(\ker M^{n_k-1})
& \dotsb\nonumber\\ & \dotsb & \supseteq
M^{n_k-1}V=\phi(U_k^{\circ}\otimes e_{n_k-1})=V_k^{\circ}\otimes
e_{n_k-1}\label{eq:(3)}.
\end{eqnarray}
Since $V$ is a semisimple $\wk$-module, taking into account
\eqref{eq:(1)}, we get from \eqref{eq:(3)}:
\bbe\label{eq:(4)}
V=(\bigoplus_{j=0}^{n_k-2}A_j)\oplus\phi(U_k),
\ee
where each $A_j$ is a complement of $M^{j+1}V$ in $M^j(\ker
M^{n_k-1})$. Analogously, taking into account \eqref{eq:(2)}, we
get from \eqref{eq:(3)}:
\bbe\label{eq:(5)}
V=(\bigoplus_{j=0}^{n_k-2}A_j)\oplus V_k.
\ee
Combining \eqref{eq:(4)} and \eqref{eq:(5)}, we see that
$\pi_k\circ\phi(U_k)=V_k$, hence $\psi_k$ is an isomorphism. To be
able to apply the inductive step, note that $A_j$'s can be chosen
in such a way that
$$
\bigoplus_{j=0}^{n_k-2}A_j=\bigoplus_{i=0}^{k-1}V_i.
$$
Indeed, this follows from the following formulas:
\begin{eqnarray*}
\ker M^{n_k-1} & = & MV\oplus (\bigoplus_{i=0}^{n_k-1}(V_i^{\circ}\otimes e_0))\quad\text{and} \nonumber \\
M^j(\ker M^{n_k-1}) & = & M^{j+1}V\oplus
(\bigoplus_{i=0}^{n_k-1}(V_i^{\circ}\otimes e_j)),\quad 0\leq
j\leq n_k-1,
\end{eqnarray*}
where each $V_i^{\circ}$ is a representation space of $\be_i$.
Thus, by \eqref{eq:(4)}
$$
V=(\bigoplus_{i=0}^{k-1}V_i)\oplus\phi(U_k)=\phi(\bigoplus_{i=0}^{k-1}U_i)\oplus\phi(U_k),
$$
which implies that the projection of
$\phi(\bigoplus_{i=0}^{k-1}U_i)$ onto $\bigoplus_{i=0}^{k-1}V_i$
is an isomorphism, hence by induction $\psi_1,\ldots,\psi_{k-1}$
are isomorphisms.
\end{proof}

As in Lemma \ref{l:sympl} above by a {\it minimal} unitary,
orthogonal, or symplectic representation we mean a unitary,
orthogonal, or symplectic representation respectively that cannot
be written as an orthogonal sum of nonzero invariant
subrepresentations.
\begin{prop}
Let $\sg'$ be a minimal unitary, orthogonal, or symplectic
admissible representation of $\wdk$. Let $U$ be a representation
space of $\sg'$ and $\langle\cdot\,,\cdot\rangle$ a nondegenerate
invariant form on $U$. Then either $\sg'$ is indecomposable or
$U\cong V\oplus\tilde V$, where $V$ is an indecomposable submodule
of $U$, $\tilde V=\check V$ if $\langle\cdot\,,\cdot\rangle$ is
bilinear, and $\tilde V=V^*$ if $\langle\cdot\,,\cdot\rangle$ is
sesquilinear. Moreover, if $\la$ is the isomorphism of
$V\oplus\tilde V$ onto $U$ and $\langle\cdot\,,\cdot\rangle '$ is
the form on $V\oplus\tilde V$ given by
$$
\langle x,y\rangle '=\langle \la(x),\la(y)\rangle,\quad x,y\in
V\oplus\tilde V,
$$
then $\langle\cdot\,,\cdot\rangle '\vert_V$ and
$\langle\cdot\,,\cdot\rangle '\vert_{\tilde V}$ are degenerate,
$\langle\cdot\,,\cdot\rangle ':V\times\tilde V\rar\bC$ is the
standard form given by
$$
\langle u,f\rangle '=f(u),\quad u\in V,\,f\in\tilde V.
$$
\end{prop}

\begin{proof}
Since $\sg'$ is minimal, it follows from Proposition \ref{pr:lok},
that $\sg'\cong\al\otimes\spin(n)$, where $\al$ is a
representation of $\wk$. Since $\sg'$ is admissible, $\al$ is
semisimple, hence $\al=\bigoplus_{i=1}^k\al_i$, where each $\al_i$
is an irreducible subrepresentation of $\al$. For each $i$ let
$U_i$ be a representation space of $\al_i\otimes\spin(n)$, so that
$U=U_1\oplus\dotsb\oplus U_k$. Let $\phi:U\rar\tilde
U_1\oplus\dotsb\oplus\tilde U_k$ be the composition of the
isomorphism induced by $\langle\cdot\,,\cdot\rangle$ with the
usual isomorphism of ${(U_1\oplus\dotsb\oplus U_k)}^{\sim}$ onto $\tilde
U_1\oplus\dotsb\oplus\tilde U_k$. For each $i$ and $j$ let
$\phi_{ij}:U_i\rar\tilde U_j$ be defined by the following diagram:
$$
\xymatrix{
U \ar[r]^-{\phi}   &   \tilde U_1\oplus\dotsb\oplus\tilde U_k \ar[d]^{\pi_j}\\
U_i \ar@{^{(}->}[u] \ar[r]^{\phi_{ij}} & \tilde U_j }
$$
where $\pi_j$ is the projection onto $j$-th factor. We claim that
for any $i$ there is $j=j(i)$ such that $\phi_{ij}$ is an
isomorphism. Indeed, let $U_i^{\circ}$ be a representation space
of $\al_i$ so that $U_i=U_i^{\circ}\otimes\bC^n$, where $\bC^n$ is
the representation space of $\spin(n)$. Let $W=\bC$ be the
representation space of $\om^{-(n-1)}$ and
$\psi:U\rar\bigoplus_{i=1}^k(\tilde U_i^{\circ}\otimes
W\otimes\bC^n)$ the composition of $\phi$ with the usual
isomorphism induced by
${\al_i\otimes\spin(n)}^{\sim}\cong\tilde\al_i\otimes\om^{-(n-1)}\otimes\spin(n)$,
$1\leq i\leq k$. For each $i$ and $j$ let
$\psi_{ij}:U_i^{\circ}\otimes\bC^n\rar\tilde U_j^{\circ}\otimes
W\otimes\bC^n$ be defined by the following diagram:
$$
\xymatrix{
U \ar[r]^-{\psi}   &   (\tilde U_1^{\circ}\otimes
W\otimes\bC^n)\oplus\dotsb\oplus (\tilde U_k^{\circ}\otimes
W\otimes\bC^n)
\ar[d]^{\pi_j}\\
U_i^{\circ}\otimes\bC^n \ar@{^{(}->}[u] \ar[r]^{\psi_{ij}} &
\tilde U_j^{\circ}\otimes W\otimes\bC^n. }
$$
Let $e_0,\ldots,e_{n-1}$ be the standard basis of $\bC^n$. 
If for
each $i$ there is $j=j(i)$ such that the projection of
$\psi(U_i^{\circ}\otimes e_0)$ onto $\tilde U_j^{\circ}\otimes
W\otimes e_0$ is nonzero, then $U_i^{\circ}\cong \tilde
U_j^{\circ}\otimes W$ (because each $U_i^{\circ}$ is irreducible),
hence $U_i\cong\tilde U_j$ and $\phi_{ij}\ne 0$. Then it follows
from Schur's lemma for indecomposable representations of $\wdk$
(\cite{r2}, p. 133, Cor. 1), that $\phi_{ij}$ is an isomorphism.

Assume now  that there is $i$ such that the projection of
$\psi(U_i^{\circ}\otimes e_0)$ onto $\tilde U_j^{\circ}\otimes
W\otimes e_0$ is zero for any $j$. Let $N$ (resp. $M$) be the
nilpotent endomorphism of $U$ (resp. of $(\tilde
U_1^{\circ}\otimes W\otimes\bC^n)\oplus\dotsb\oplus (\tilde
U_k^{\circ}\otimes W\otimes\bC^n)$). Then $\psi(U_i^{\circ}\otimes
e_0)\subseteq X$, where $X=\bigoplus_{t\geq 1;s}(\tilde
U_s^{\circ}\otimes W\otimes e_t)$ and $X\subseteq\ker M^{n-1}$.
Since $U_i^{\circ}\otimes e_0\not\subseteq\ker N^{n-1}$, we get a
contradiction with $\psi$ being an isomorphism.

Thus, in particular, there exists some $j$ such that $\phi_{1j}$
is an isomorphism. If $j=1$ then
$\langle\cdot\,,\cdot\rangle\vert_{U_1}$ is nondegenerate, hence
$U_1$ and its orthogonal complement are invariant subspaces of
$U$. Since $U$ is minimal, it implies that $U=U_1$ and $U$ is
indecomposable. If $j\ne 1$ then without loss of generality we can
assume that $j=2$, $\langle\cdot\,,\cdot\rangle\vert_{U_1}$ and
$\langle\cdot\,,\cdot\rangle\vert_{U_2}$ are degenerate. Let us
show that $\langle\cdot\,,\cdot\rangle\vert_{U_1\oplus U_2}$ is
nondegenerate. Suppose it is degenerate, i.e.
$K=\ker(\langle\cdot\,,\cdot\rangle\vert_{U_1\oplus U_2})$ is
nonzero. Let $R_1$ (resp. $R_2$) be the nilpotent endomorphism of
$U_1$ (resp. $U_2$). Then $R=R_1\oplus R_2$ is the nilpotent
endomorphism of $U_1\oplus U_2$. We claim that $K\bigcap\ker R\ne
0$. Indeed, let $x\in K$ and $x\ne 0$. Then there is $i$ $(0\leq
i\leq n-1)$ such that $R^ix\in\ker R$ and $R^ix\ne 0$. Also,
$$
\langle R^ix,y\rangle=(-1)^i\cdot\langle
x,R^iy\rangle=0\quad\text{for any }y\in U_1\oplus U_2,
$$
hence $R^ix\in K$. Let $x\in K\bigcap\ker R$ and $x\ne 0$, i.e.
$x=x_1+x_2$, where $x_i\in\ker R_i$, $i=1,2$. Without loss of
generality we can assume that $x_1\ne 0$. Since $\phi_{12}$ is an
isomorphism, there is $y_2\in U_2$ such that $\langle
x_1,y_2\rangle\ne 0$. By assumption,
$\langle\cdot\,,\cdot\rangle\vert_{U_2}$ is degenerate, hence
$K_2=\ker(\langle\cdot\,,\cdot\rangle\vert_{U_2})$ is nonzero.
Then by the same argument as above $K_2\bigcap\ker R_2\ne 0$. Since $\ker
R_2=U_2^{\circ}\otimes e_{n-1}$, it is irreducible, consequently
$\ker R_2\subseteq K_2$. In  particular, $\langle
x_2,y_2\rangle=0$. Hence
$$
\langle x_1+x_2,y_2\rangle=\langle x_1,y_2\rangle.
$$
Since $\langle x_1,y_2\rangle\ne 0$ by the choice of $y_2$, we get
a contradiction with $x_1+x_2\in K$. Thus,
$\langle\cdot\,,\cdot\rangle\vert_{U_1\oplus U_2}$ is
nondegenerate. Since $U$ is minimal the same argument as above
implies that $U=U_1\oplus U_2\cong \tilde U_2\oplus U_2$.
\end{proof}

\end{document}